\documentclass[conference]{IEEEtran}
\usepackage{times}

\makeatletter
\def\ps@headings{%
\def\@oddhead{\mbox{}\scriptsize\rightmark \hfil \thepage}%
\def\@evenhead{\scriptsize\thepage \hfil \leftmark\mbox{}}%
\def\@oddfoot{}%
\def\@evenfoot{}}
\makeatother
\usepackage{amsfonts,amsmath,amstext,verbatim}
\usepackage{amssymb}%
\usepackage{setspace}
\usepackage{ifpdf}

\makeatother  \renewenvironment{abstract}{%
  \small\bfseries\textit{Abstract}:  }

\ifCLASSINFOpdf
\usepackage[pdftex]{graphicx}
\usepackage[pdftex,draft,
          bookmarks=false,%
          bookmarksnumbered=false,%
          pdftitle=Connectivity\ in\ Sub-Poisson\ Networks,%
          pdfauthor=B.\ Blaszczyszyn%
          \ -\  D.\ Yogeshwaran]{hyperref} 
\usepackage[numbers,sort&compress]{natbib}
\usepackage{hypernat}
\graphicspath{{./figures/}{./}}
\DeclareGraphicsExtensions{.pdf}
\else

\usepackage[dvips]{graphicx}
\graphicspath{{/.}{./figures/}}
\DeclareGraphicsExtensions{.eps}
\usepackage[numbers,sort&compress]{natbib}
\fi

\newcommand{\md}{\text{\rm d}}
\newcommand{\ir}{{\mathbb{R}}}
\newcommand{\mR}{{\mathbb{R}}}
\newcommand{\E}{{\mathbf E}}
\newcommand{\sE}{{\mathbf E}}
\newcommand{\Pro}{{\mathbf P}}
\newcommand{\ind}{1\hspace{-0.30em}{\mbox{I}}}

\newcommand{\lam}{\lambda}

\renewcommand{\theequation}{\arabic{section}.\arabic{equation}}

\newtheorem{Th}{Theorem}[section]
\newtheorem{prop}[Th]{Proposition}
\newenvironment{Prop}{\bf\begin{prop}\rm\em}{\end{prop}} 
\newtheorem{res}[Th]{Result}
\newtheorem{fact}[Th]{Fact}
\newtheorem{ex}[Th]{Example}

\newcommand{\calO}{{\mathcal{O}}}

\newcommand{\bbZ}{\Bbb{Z}}
\newcommand{\bbR}{\Bbb{R}}

\begin{document}
\title{Connectivity in Sub-Poisson Networks\vspace{-1ex}}

\author{
\large {\bf Bart{\l }omiej B{\l }aszczyszyn}
\ \ and\ \  {\bf D. Yogeshwaran}\\
INRIA-ENS, 23 Avenue d'Italie  75214 Paris, France\\
\{Bartek.Blaszczyszyn,  Yogeshwaran.Dhandapani\}@ens.fr
\vspace{-2ex}
}

\maketitle

\thispagestyle{empty} \pagestyle{empty}

\begin{abstract}
We consider a class  of point processes (pp), which we call {\em
  sub-Poisson}; these are pp
that can be  directionally-convexly ($dcx$) dominated 
by some Poisson pp. 
The $dcx$ order has already been shown in~\cite{snorder} useful 
in comparing various point process characteristics,
including Ripley's and correlation functions 
as well as shot-noise fields generated by pp,
indicating in particular 
that smaller in the $dcx$ order 
processes exhibit more regularity (less
clustering, less voids) in the repartition of their points. 
Using these results, 
in this paper we study the impact of the $dcx$ ordering of pp
on the properties of two  continuum percolation models, 
which have been proposed in the literature
to address macroscopic connectivity properties of
large wireless networks.
As the first main result of this paper, 
we extend the classical result on the existence
of phase transition in the percolation of the Gilbert's graph
(called also the Boolean model), generated by a homogeneous Poisson
pp,  to the class of homogeneous sub-Poisson pp.
We also extend  a recent result of the same nature for the
SINR graph, to sub-Poisson pp.
Finally, as examples we show that the so-called 
perturbed lattices are sub-Poisson.  
More generally,  perturbed lattices  provide some spectrum of models 
that ranges from periodic grids, usually considered in cellular network
context, to Poisson ad-hoc networks, 
and to various more clustered pp including
some doubly stochastic Poisson ones. 
\end{abstract}

\begin{keywords} 
percolation, $dcx$ order,  Gilbert's graph, Boolean model, SINR graph, 
 wireless network, Poisson point process, perturbed lattice, 
determinantal point process, connectivity, capacity  
\end{keywords}

\vspace{-2ex}\section{Introduction}
A {\em network}, in the simplest terms, 
is a collection of points in some space (e.g. on the Euclidean plane),
called {\em nodes} or {\em vertexes}, 
and a collection of node pairs, called {\em edges}. 
The presence of an edge between two nodes indicates that they can
directly communicate with each other.
The mathematical name for this network model is {\em graph}.

A class of networks that has
recently attracted particular interest in the wireless communication context, 
is called {\em ad-hoc networks}. It is distinguished  by the fact that
the network  nodes  are not subject to any regular (say periodic) 
geometric emplacement in the space
but can be rather seen as a snapshot
of some {\em random point pattern}, called  also {\em point process}
(pp)  in
the mathematical formalism typically used in this context.

Connectivity, i.e.  possibility of indirect, multi-hop,
communication between distant  
nodes,  is probably the first issue that has to be addressed
when considering ad-hoc networks.
An ubiquitous assumption when studying this problem 
is that the node  randomness is modeled by a
spatial  {\em Poisson} pp. This latter situation can be
characterized by independence and Poisson distribution of the number
of nodes observed in disjoint subsets of the space.  
Poisson assumption  in the above context is often too simplistic,
however analysis or even modeling of networks without this assumption 
is in most cases very difficult.

In this paper we  introduce some class of point
processes that might be roughly described as exhibiting 
less variable point patterns than Poisson pp. 
We call point processes (pp's) 
of this class {\em sub-Poisson point processes}. 
Our objective is two-fold. On one hand we want to argue that this class of pp's allows to 
extend some classical results regarding network connectivity
and on the other hand we want to bring  attention to sub-Poisson pp's, 
as they  might be useful for modeling of ad-hoc networks, 
whose nodes are more regularly distributed than a Poisson~pp.

\subsection{Sub-Poisson Random Variables}
A random variable is called sub-Poisson, 
if its variance is not larger  than its mean (with the equality
holding true for Poisson variable).
Intuitively, if strictly sub-Poisson variables
were to describe the number of nodes in different subsets of the space
then the resulting point patterns would exhibit less clustering (or
bunching) than the Poisson point pattern having the same number of
points per unit of space volume. This in turn, still intuitively,
should have positive impact on the connectivity and perhaps  capacity
and  other network performance metrics, the reason being that
the perfectly regular, periodic patterns are commonly considered (and
sometimes can be proved) as being optimal. The aim of the present
article is to provide rigorous results on the comparison of the
connectivity and capacity properties of certain 
sub-Poisson networks to these  of the respective Poisson networks.

\subsection{$dcx$ Sub-Poisson  Point Processes}
The statistical variability of random variables (say with the same
mean) can be compared   
only to some limited extent by looking at their  variances, 
but more fully by convex ordering.
Under this order one can compare the expected values of {\em all} convex
functions of these variables. 
In multi-dimensions there is no one single notion of convexity.
Besides different statistical
variability of marginal distributions, two random vectors
(think of number of nodes in different subsets of the space) can
exhibit different dependence properties on their coordinates.
The most evident example here is comparison of the  vector
composed of several copies of one random variable to a vector
composed of independent copies sampled from the same distribution.
Among several notions of ``convex-like'' ordering of random vectors
the so called {\em directionally convex} ($dcx$) order, allowing one  to compare
expectations of all $dcx$ functions (see Section~\ref{sec:dcx_defns}
below) of these vectors,
takes into account both the  dependence structure of random
vectors and  the variability of their marginals.
It can be naturally extended to random fields by
comparison of all finite dimensional distributions,
as well as to random pp'as and even locally finite random 
measures by  viewing them as non-negative fields of
measure-values on all bounded Borel subsets of the space;
cf.~\cite{snorder}. 

Using this latter formalism, we say that a  pp
is {\em \hbox{$dcx$-}sub-Poisson}  (or simply {\em sub-Poisson} when there
is no ambiguity) if it is $dcx$-smaller than a Poisson
pp having the same mean measure, i.e., 
the mean number of nodes in any given set.
We will also say {\em sub-Poisson network} in the case when 
the nodes of the graph modeling  the network are distributed according
to some sub-Poisson pp.

We shall also  see that there are classes of pp whose so-called 
{\em joint intensities} (densities of the higher order moment
measures, when they 
exist;  see Section~\ref{ss.intensities})  are
smaller than those of Poisson pp. Using this latter (and weaker) property
we define the class of {\em weakly sub-Poisson point processes}.  


Another weakening of the comparability assumption, to the expectations of
increasing (or decreasing) $dcx$ functions 
allows to define  {\em increasing $dcx$} sub-Poisson pp's, which 
can be seen as being less variable than some Poisson pp 
having  possibly smaller (or larger) mean measure. We shall abbreviate increasing $dcx$ by $idcx$ and 
similarly {\em decreasing} $dcx$ by $ddcx$.

Our choice of the $dcx$ order to define sub-Poisson  pp's 
has its roots in~\cite{snorder}, where one shows various results as
well as examples indicating that the $dcx$ order implies ordering
of several  clustering characteristics known in spatial statistics
such as Ripley's K-function or second moment densities.
Namely, a  pp that is larger in the $dcx$ order exhibits
more clustering (while having the same mean number of points in any given set).

\subsection{Connectivity of Sub-Poisson Networks}

Full connectivity  (multi-hop communication between any two nodes)
of an ad-hoc network with many nodes
is typically hard to maintain, and so, a more modest
question of existence of a large enough, connected subset of nodes
(called component) is studied. A possible approach to this  problem,
proposed in~\cite{Gilbert61}, and based on the mathematical 
theory of {\em percolation}, consists in studying existence of an
{\em infinite} (called {\em giant}) {\em component} of the  
infinite graph modeling a network.
Existence of such a component 
is interpreted as an indication that the connectivity of 
the modeled  ad-hoc network scales well with its number of nodes.

\subsubsection{Percolation of Gilbert's  Network}
\label{sss.Gilbert}
Percolation models have been extensively studied both in mathematical
and communication literature.
The model proposed in~\cite{Gilbert61},
called now the {\em Gilbert's model}~\footnote{also {\em Boolean
    model} or {\em random geometric graph}}, 
is now considered as the classical continuum model in
percolation theory. 
It assumes that each node 
has a given fixed range of communication $\rho$ and direct connection
between any two nodes is feasible if they are within the distance $\rho$ 
from  each other, regardless of the positions of other nodes in the network. 
The known answer  to the percolation problem in this model is  given  under
the assumption that the nodes are distributed according to a homogeneous 
Poisson pp having a {\em density} of $\lambda$ nodes on
average per unit of volume (or surface). The result   says that
there exists a {\em non-degenerate critical communication range}
$0<\rho_c=\rho_c(\lambda)<\infty$,
such that for $\rho\le \rho_c$ there is no giant connected
component of the network (the model does not  percolate), 
while for $\rho>\rho_c$ there is exactly one such
component (the model percolates), both statements holding true almost
surely; i.e., for 
almost all Poisson realizations of the network.
An equivalent statement of the above result says that for a given
communication  range $\rho$ there exists a non-degenerate
critical density of nodes
$0<\lambda_c=\lambda_c(\rho)<\infty$, below which the model does not
percolate and above which it does so almost surely. 

As one of the main results of this paper we will prove an extension of
the above result, which says that the {\em critical communication range of
any homogeneous sub-Poisson Gilbert's network model
is not degenerate}. Moreover,
this critical communication range is
bounded away from zero and infinity by the constants which depend only
on the mean density of nodes and not on the finer structure of the
sub-Poisson pp of network nodes. 
 Partial results, regarding only non-degeneracy
at zero or at infinity can be proved for $idcx-$ and
$ddcx$-sub Poisson networks respectively.

\subsubsection{Percolation of SINR Networks}
A more adequate percolation model of a wireless communication network,
called the {\em SINR} graph,
was studied more recently in~\cite{Dousse_etal_TON,Dousse_etal}. It allows 
one to take into account the interference intrinsically related to
wireless communications.
The interference power is modeled by the {\em
shot-noise} field generated by the pp 
of transmitting nodes. Each pair of nodes 
in the new model is joined by an edge when  the signal
power to this shot-noise plus some other (external) noise power ratio 
is large enough. 
The resulting  random graph does not have the independence structure 
of the Gilbert's  model and increasing the communication range
(equivalent in this model to increasing the signal power and
hence the value of the shot-noise) is not necessarily 
beneficial for connectivity. Similar observation holds as regards the
increase in the node density. 
In fact, the above SINR network model has two essential parameters: 
the density   $\lambda$ of Poisson pp of nodes and the   
shot-noise reduction factor $\gamma$. The percolation domain is characterized 
in the Cartesian product of these two parameters. 
The key result of \cite{Dousse_etal} says that whenever $\lambda$ is larger than the 
critical value corresponding to the percolation of the model with interference perfectly
canceled out (for $\gamma=0$, which simplifies the model to
the Gilbert's one) then there exists a critical value
$\gamma_c=\gamma_c(\lambda)>0$ such that the model
percolates for $\gamma<\gamma_c$ and does not percolate for
$\gamma>\gamma_c$. As the second main result of this paper, we  extend
the above result, also  to sub-Poisson SINR networks.

\subsection{Comparison of Shot-Noise Fields}
Our proofs rely on the comparison of extremal and additive shot-noise
fields generated by $dcx$-ordered pp's studied
in~\cite{snorder}. While the connection to the SINR model is evident
(additive shot-noise is an element of this model) the 
connection to the Gilbert's model is perhaps less evident and relies
on the fact that this latter model  can be represented as a upper level-set 
of some extremal shot-noise field.

More precisely, from~\cite[Propsition~4.1]{snorder}
one can conclude  that the probability of $n$ given locations 
in the space {\em not} being within the communication
range of any of the nodes in the Gilbert's model 
is higher for the network whose nodes are modeled by a pp
larger in $dcx$ order. Using this property, suitable discretization of the model
and the Peierls argument 
(cf.~\cite[ pp.~17--18]{Gr99} or~\cite[Proposition 14.1.4]{FnT1}) 
one can prove finiteness of the critical transmission range in the
sub-Poisson Gilbert's network. The strict positivity of this range 
can be proved by comparing the expected number of paths from the origin
to the boundary of an increasing box, again  in some suitable
discretization of the model. This latter comparison can be done relying only on weak
sub-Poisson assumption.

{\em The remaining part of this paper is organized as follows}: In
Section~\ref{s.Preliminaries} we provide necessary notions and
notation. 
Sub-Poisson Gilbert's model  is studied in Section~\ref{s.Gilberts}
and sub-Poisson SINR model in Section~\ref{s.SINR}. 
In Section~\ref{s.Capacity} we make some remarks on 
the impact of $dcx$ ordering on various model characteristics usually called 
capacities, including  the so-called capacity functional, being one of the
fundamental characteristics studied in stochastic geometry, as well as
on some other capacity quantifiers in the information-theoretic sense.
In Section~\ref{s.Examples} we show some examples of sub-Poisson and
weakly sub-Poisson pp, in particular the so-called {\em perturbed
lattices}.  Conclusions as well
as open questions are presented in Section~\ref{s.Conclusions}.  

\section{Preliminaries}
\label{s.Preliminaries}

\subsection{$dcx$ Order}
\label{sec:dcx_defns}
We say that a function $f:\mR^d \to \mR$ is
{\em directionally convex}~($dcx$) if for every $x,y,p,q \in \mR^d$ such that
$p \leq x,y \leq q$ (i.e., $p \leq x \leq q$ and $p \leq y \leq q$ with inequalities understood component-wise)
and $x+y = p+q$ one has  $f(x) + f(y) \leq f(p) + f(q)$. 
Also, we shall abbreviate {\em increasing} and $dcx$ functions by $idcx$ and
{\em decreasing} and $dcx$ by $ddcx$. 

For two  real-valued random vectors $X$ and $Y$  of the
same dimension,   $X$ is said to be $dcx$  smaller than $Y$,
(denoted by $X\le_{dcx}Y$) 
if $\sE(f(X)) \leq \sE(f(Y))$ for every $dcx$ function $f$
for which both expectations are finite. 
In full analogy one defines $idcx$ and $ddcx$ orders of random vectors
considering $idcx$ and $ddcx$ functions, respectively.
These orders clearly depend only on the distributions of the vectors. 
Two real valued stochastic processes (or random fields) 
are said $dcx, idcx$ or $ddcx$ ordered 
if any finite-dimensional distributions of these processes are ordered.
This definition extends also to locally finite random measures (in particular
pp's), by ordering random values of these measures (in particular
numbers of points) on any finite collection of subsets of the state space. More precisely, we say for two pp's $\Phi_1,\Phi_2$ that $\Phi_1 \leq_{dcx(idcx,ddcx)} \Phi_2$ if for every finitely many bounded Borel subsets $B_1,\ldots,B_n$ we have that,
$$(\Phi_1(B_1),\ldots,\Phi_1(B_n)) \leq_{dcx(idcx,ddcx)} (\Phi_2(B_1),\ldots,\Phi_2(B_n))  .$$
It was shown in~\cite{snorder} that  verifying  the above
property for all  {\em mutually disjoint bounded Borel subsets} $B_i$ is 
a sufficient condition for the respective ordering of the pp's.

Using the above definition, we say that  {\em $\Phi$ is sub-Poisson if $\Phi
\leq_{dcx} \Phi_{Poi}$, where $\Phi_{Poi}$ is some Poisson pp}.
Noting that for the  linear function $f(x)=x$, both $f$ and
$-f$ is convex (and thus $dcx$ in one dimension) 
one can observe that $\Phi_{Poi}$ needs to have the same mean measure 
as $\Phi$; i.e., if  $\Phi
\leq_{dcx} \Phi_{Poi}$  then $\E[\Phi(B)]=\E[\Phi_{Poi}(B)]$ for every
bounded Borel subset~$B$. In particular, {\em a stationary pp $\Phi$ is
sub-Poisson  if it is $dcx$ smaller than Poisson pp $\Phi_\lambda$ 
of the same intensity} $\lambda=\E[\Phi(B)]/|B|$, where $|B|$ is the
$d\;-$dimensional Lebesgue measure of the set $B$.   
{From} now on $\Phi_\lambda$ or $\Phi_\mu$ will always denote
homogeneous Poisson pp of intensity $\lambda$ or $\mu$.
We will also say that pp $\Phi$ is homogeneous if its mean measure is 
equal, up to a constant, to the Lebesque measure; i.e. 
$\E[\Phi(B)]=\lambda |B|$, for some constant $\lambda$
all bounded Borel sets $B$. This is a
weaker assumption than the stationarity of $\Phi$.
 
We say that $\Phi$ is $idcx$($ddcx$)-sub-Poisson 
if $\Phi\leq_{idcx(ddcx)}\Phi_{Poi}$,  where $\Phi_{Poi}$ is some
Poisson pp. In this case the mean measure of $\Phi$ is smaller or
equal (larger or equal)  to that of $\Phi_{Poi}$.

\subsection{Joint Intensities}
\label{ss.intensities}
The joint intensities $\rho^{(k)}(x_1,\ldots,x_k)$ of a pp are defined by the following relation for every finitely many disjoint bounded Borel sets $B_1,\ldots,B_n$:
$$ \sE\Bigl[\prod_{i=1}^k\Phi(B_i)\Bigr] = 
\int\!\!\dots\!\!\int_{B_1\times\dots\times B_k}
\rho^{(k)}(x_1,\ldots,x_k) \md x_1 \ldots \md x_k,$$
provided  they exist (i.e., the respective moment measures admit
densities). 
 A pp $\Phi$ with joint intensities of all orders $k$ is said to be
 {\em weakly sub-Poisson} if there exists a constant $\lam$ such that
 for all $k \geq 1$,  
\begin{equation}\label{e.weakly-sub-Poi}
 \rho^{(k)}(x_1,\ldots,x_k) \leq \lam^k \, \, \, \, \,\mbox{a.e.}.
\end{equation}
Due to the $idcx$ property of the function $f(x_1,\ldots,x_k)=
\prod_i x_i^+$, the ordering  $\Phi_1
\leq_{idcx} \Phi_2$ implies that for all $k \geq 1$,
$\sE[\prod_{i=1}^k\Phi_1(B_i)] \leq \sE[\prod_{i=1}^k\Phi_2(B_i)]$ for
any disjoint bounded Borel sets $B_1,\ldots,B_k$. Hence, the respective joint
intensities $\rho_j^{(k)}$ of pp's $\Phi_j$, $j=1,2$, provided they exist,
 obey for all $k \geq 1$,  $\rho^{(k)}_1(x_1,\ldots,x_k)
\leq \rho^{(k)}_2(x_1,\ldots,x_k)$ a.e.. Since 
$\rho^{(k)}(x_1,\ldots,x_k) = \lam^k$ for Poisson pp $\Phi_{\lam}$, we are
justified in using the term weakly sub-Poisson for pp $\Phi$
satisfying~(\ref{e.weakly-sub-Poi}).

\section{Percolation of Sub-Poisson Gilbert's Networks}
\label{s.Gilberts}

Given a pp $\Phi$ on $\ir^d$ and a non-negative constant
$r>0$ one defines the ({\em spherical}) Boolean model generated
by~$\Phi$ of ball of radius $r$ as the union
$C(\Phi,r)=\bigcup_{X\in\Phi}B_{X}(r)$, where $B_x(r)$ is the ball
centered at $x\in\ir^d$ of radius~$r$. We say that $C(\Phi,r)$
{\em percolates} if there exists an unbounded, connected subset of
$C(\Phi,r)$. This definition extends to any random subset of $\ir^d$.
One defines the {\em critical radius} $r_c(\Phi)$ of the Boolean model 
as $r_c(\Phi):=\inf\{r>0: \Pro\{\, C(\Phi,r)\;\text{percolates}\}>0\}$. 
Here is the first result of this paper. 
\begin{Prop}
\label{thm:perc_sub-Poisson_pp} 
There exist universal constants $0<c(\lambda)$ and $C(\lambda)<\infty$
depending only on $\lambda$ (and the dimension $d$) 
such that if $\Phi \leq_{idcx} \Phi_{\lambda}$ then $c(\lambda) \leq
r_c(\Phi)$ and if $\Phi \leq_{ddcx} \Phi_{\lambda}$ then 
$r_c(\Phi) \leq C(\lambda)$.
Thus, for a homogeneous $dcx$-sub-Poisson $\Phi$ of intensity $\lambda$ we have
$ 0 < c(\lambda) \leq r_c(\Phi)  \leq C(\lambda) <
\infty $, where the constants depend only on $\lambda$.
The lower bound holds  also when $\Phi$ is weakly sub-Poisson.
\end{Prop}
The proof is given in the Appendix.
The above result  can be extended to the so-called  $k$-percolation models. 

Note that the connectivity structure of the Boolean model $C(\Phi,r)$
corresponds to that  of the Gilbert's network with nodes in
$\Phi$ and communication range $\rho=2r$;
cf. Section~\ref{sss.Gilbert}. Thus the critical
communication range  $\rho_c=\rho_c(\Phi)$ in this latter network
is twice the critical
radius of the corresponding Boolean model.

While the finiteness of the critical radius of the Boolean model
(and thus communication range of the sub-Poisson Gilbert's network) is
intuitively a desired property, its positivity 
at first glance might be seen from the networking point of view
as irrelevant, if not a disadvantage.
A deeper inspection of wireless communication mechanisms  shows however
that sometimes  a non-percolation might be also a desired property. 
The following modification of the Gilbert's model, 
that can be seen as a toy version of the SINR model studied in
Section~\ref{s.SINR}, sheds  some light on this latter statement.

\begin{ex}[Gilbert's carrier-sense network]
Consider a planar 
ad-hoc network consisting of nodes modeled by a point
process $\Phi_B$ on the plane (in $\ir^2$) 
and  having communication range $\rho$.
This process corresponds to a {\em back-bone} of the  network,
whose percolation we are looking for.
Consider also an auxiliary pp $\Phi_I$ of {\em interferers}
also on $\ir^2$. Consider the following modification of the Gilbert's
connection rule: any two given nodes of $\Phi_B$ can directly communicate
(are joined be an edge) when they are within the communication range
$\rho$ from each other,  
however only when there is no interfering node (point  of $\Phi_I$)
within the {\em  sensing range} $R>\rho$ of any of these two nodes.
Note that any connected component of this network is included in some
connected subset of the {\em complement} $\calO=\ir^d\setminus C(\Phi_I,R-\rho)$
of the Boolean model $C(\Phi_I,R-\rho)$.  Percolation of
this {\em vacant region} $\calO$ of the spherical Boolean model
is thus a necessary condition for the percolation of our  modification 
of the Gilbert's network. Now, the vacant region of the planar Boolean model 
cannot percolate if the Boolean model itself does (cf~\cite[Theorem
  4.4]{MR96}). Thus a non-percolation of some Boolean model 
related to the interferes $\Phi_I$ is a necessary condition for 
the percolation of the communication network on $\Phi_B$.
\end{ex}

\section{Percolation of Sub-Poisson SINR Networks}
\label{s.SINR}

In this section, we shall work only on the plane $\mR^2$.
We slightly modify the definition of SINR network introduced
in~\cite{Dousse_etal_TON} allowing for {\em external interferers}.
The parameters of the model are non-negative numbers $P$ (signal
power),  $N$ (environmental noise), $\gamma$ (interference reduction
factor), $T$ (SINR threshold) and an attenuation 
function $l(r)$ of the distance $r\ge0$ satisfying
$0\le l(r)\leq 1$, continuous,  strictly decreasing~\footnote{So
  it is rather path {\em gain} function.}
 on its support, with 
$l(0) \geq \frac{TN}{P}$ and 
$\int_0^{\infty} x l(x) \md x  < \infty$.
These assumptions are exactly as in~\cite{Dousse_etal}.

Given a pp $\Phi$, the (unit-power) interference generated by $\Phi$
at location $x$ is defined as 
$I_{\Phi}(x) := \sum_{X \in \Phi\setminus\{x\}} l(|X-x|)$.
More generally this object is also  called  (additive) shot-noise field 
generated on by $\Phi$ withe {\em response  function}
$l$; cf~\cite[Ch.~2.2 and~2.3]{FnT1}. 
The SINR from $x$ to $y$ with interference from $\Phi$ is defined as 
\begin{equation}
\label{eqn:sinr_defn}
\text{SINR}(x,y,\Phi,\gamma) :=  \frac{Pl(|x-y|)}{N + \gamma P I_{\Phi\setminus\{x\}}(y)}.
\end{equation}
Let $\Phi_B$ and $\Phi_I$ be two pp's on the plane. 
We do not assume any particular dependence 
between  $\Phi_B$ and $\Phi_I$. In particular one may think of
$\Phi_B\subset\Phi_I$. 
Let $P,N,T > 0$ and $\gamma \geq 0$. 
The SINR network with {\em back-bone}
$\Phi_B$ and {\em interferers} $\Phi_I$ is defined as
a graph $G(\Phi_B,\Phi_I,\gamma)$ with nodes in $\Phi_B$ and 
edges joining any two nodes $X,Y\in\Phi_B$ when 
$\text{SINR}(Y,X,\Phi_I,\gamma) > T$ and
$\text{SINR}(X,Y,\Phi_I,\gamma) > T$.
The  SNR graph (i.e, the graph without interference)
is defined as $G(\Phi_B)=G(\Phi_B,\Phi_I,0)$, which is equivalent also 
to taking $\Phi_I = \emptyset$.
Observe that the SNR graph $G(\Phi_B)$ corresponds to the Gilbert's
network with nodes in $\Phi_B$ of communication range 
$\rho_l = l^{-1}(\frac{TN}{P})$. Percolation in the above graphs 
is existence of an infinite connected component in the graph-theoretic sense.

\subsubsection{Poisson Back-Bone}
\label{sec:poisson_sinr}
Firstly, we consider the case when the backbone nodes are distributed
according to  Poisson pp $\Phi_B = \Phi_{\lambda}$, for some $\lambda>0$.
We shall use $G(\lambda,\Phi_I,\gamma)$ and $G(\lambda)$ to denote the
corresponding SINR and SNR graphs respectively. 
Recall from~\ref{sss.Gilbert}
that $\lam_c(\rho)$ is the critical intensity for percolation of
the Poisson Gilbert's network of communication range~$\rho$.
The following result guarantees the existence of $\gamma > 0$ such
that for any homogeneous sub-Poisson pp of interferers 
$\Phi_I$ the SINR network  $G(\lam,\Phi,\gamma)$
percolates provided $G(\lam)$ percolates.

\begin{Prop}
\label{thm:sinr_poisson_perc} 
Let $\lam > \lam_c(\rho_l)$ and $\Phi \leq_{idcx} \Phi_{\mu}$ for some
$\mu > 0$. Then there exists $\gamma_c=\gamma_c(\lambda,\mu,P,T,N) >
0$ 
such that
$G(\lam,\Phi,\gamma)$ percolates for $\gamma<\gamma_c$. 
\end{Prop}
The proof is given in the Appendix.

Recall that we have not assumed the independence of $\Phi_I$ and
$\Phi_B=\Phi_{\lam}$. In particular, one can take $\Phi_I=\Phi_{\lam}
\cup \Phi'$ where $\Phi'$ are some {\em external} interferers. 
If $\Phi'$ is $idcx$-sub-Poisson and independent of $\Phi_\lambda$
then $\Phi_I$ is also $idcx$-sub-Poisson
(cf~\cite[Proposition~3.2]{snorder}). 
The result of Proposition~\ref{thm:sinr_poisson_perc} in the 
special case of $\Phi' = \emptyset$ was proved in \cite{Dousse_etal}. 
The present extension allows for any $idcx$-sub-Poisson 
pattern of independent external interferers. The proof of our 
result  can be also modified (which will not be presented 
in this version of the paper) to allow for external interferer's
$\Phi'$ possibly dependent of the backbone. In this full generality
the result says that {\em any   
external pattern of homogeneous $idcx$-sub-Poisson interferers added to the SINR
network of \cite{Dousse_etal} cannot make the giant component of this
network to disappear, provided the interference cancellation factor
$\gamma$ is appropriately adjusted}.  Since an $idcx$-sub-Poisson point
process can be of arbitrarily large intensity, the above observation
can be loosely rephrased in the following form: {\em It is not the density
of interferers  that matters for the network connectivity, but their
structure;  $idcx$-sub-Poisson interferers do not hurt essentially
the network connectivity}.
 
\subsubsection{Sub-Poisson Back-Bone}
\label{sec:non_poisson_sinr}

We shall now consider the case when the backbone nodes are formed by a 
sub-Poisson pp. In this case, we can give a weaker result, 
namely that with appropriately increased signal power $P$, 
the SINR graph will percolate for small interference parameter $\gamma > 0$. 
This corresponds to an early version of the result for the Poisson SINR
network, proved in~\cite{Dousse_etal_TON}, 
 where the percolation of the SINR network is guaranteed for the 
intensity of nodes possibly larger than the critical one in 
the corresponding SNR network.
\vspace{-1ex}
\begin{Prop}
\label{thm:perc_sinr_sub-Poisson} 
Let $\Phi_B \leq_{ddcx} \Phi_{\lam}$ for some $\lam > 0$ and $\Phi_I
\leq_{idcx} \Phi_{\mu}$ for some $\mu > 0$ and also assume that $l(x)
> 0$ for all $x \geq 0$. Then there exist $P, \gamma > 0$ such
that $G(\Phi_B,\Phi_I,\gamma)$ percolates. 
\end{Prop}
The proof is given in the Appendix.

As in Proposition~\ref{thm:sinr_poisson_perc}, we have not assumed the
independence of $\Phi_I$ and $\Phi_B$. In particular  one can take
$\Phi_I = \Phi_B \cup \Phi'$,
where $\Phi'$ and $\Phi_B$ are independent and $dcx$-sub-Poisson. 
%

%

\section{$dcx$  Ordering and Capacity}
\label{s.Capacity}
Now we want to make some remarks  on the relation between $dcx$
ordering and some model characteristics 
usually called  capacities, including  the so-called
stochastic-geometric capacity
functional, as well as
some other capacity quantifiers in the information-theoretic  sense.
Our main tool is the following result 
proved in~\cite[Theorem 2.1]{snorder}. It  says that
 $\Phi_1\leq_{dcx}\Phi_2$
implies the same ordering of the respective shot-noise fields 
(with an arbitrary non-negative response function) 
\begin{equation}\label{e.shot-noise-dcx}
(I_{\Phi_1}(x_1),\ldots, I_{\Phi_1}(x_n))\le_{dcx}
(I_{\Phi_2}(x_1),\ldots, I_{\Phi_2}(x_n))\,
\end{equation}
for any $x_i, i=1,\ldots,n$.
This implies in particular that 
\begin{equation}
\label{e.shot-noise-LT}
\E\Bigl[\exp\Bigl\{s\sum_{i=1}^nI_{\Phi_1}(x_i)\Bigr\}\Bigr]\le
\E\Bigl[\exp\Bigl\{s\sum_{i=1}^nI_{\Phi_2}(x_i)\Bigr\}\Bigr]
\end{equation}
for {\em both positive and negative $s$}.

\subsection{$dcx$ Ordering and  the Capacity Functional}
\label{ss.SGcapacity}
{\em Capacity functional} $T_\Xi(B)$ of a random set $\Xi$
is defined as $T_{\Xi}(B)=\Pro\{\,\Xi\cap B\not=\emptyset\,\}$ for all  bounded
Borel  sets~$B$. A fundamental result of stochastic geometry, called 
the Choquet's theorem (cf~\cite{Matheron75})  
says that the capacity functional defines entirely the
distribution of a random closed set. The complement of it,
$V_\Xi(B)=1-T_\Xi(B)$, is called {\em void probability functional}.

Inequality~(\ref{e.shot-noise-LT}) allows to compare capacity
functionals and  void probabilities of Boolean models 
generated by $dcx$ ordered pp's.
Indeed, note that 
$$
V_{C(\Phi,r)}(B)=\E\Bigl[\prod_{X\in\Phi}\!\ind(|X-B|>r)\Bigr]
=\E\Bigl[\exp\Bigl\{\sum_{X\in\Phi}h(X)\Bigr\}\!\Bigr],
$$
where $h(x)=\log(\ind(|X-B|>r))$. Note that the latter expression has
a form of the shot-noise variable and thus
using the inequality analogous to~(\ref{e.shot-noise-LT}) for shot-noise 
with the response function $h(\cdot)$ (cf~\cite[Theorem 2.1]{snorder}
for such a generalization) we observe that
if $\Phi_1\leq_{dcx}\Phi_2$ then 
\begin{eqnarray*}
V_{C(\Phi_1,r)}(B)&\le& V_{C(\Phi_2,r)}(B)\\
T_{C(\Phi_1,r)}(B)&\ge& T_{C(\Phi_2,r)}(B)\,.
\end{eqnarray*}
In other words, one can say that $dcx$ smaller pp's exhibit smaller
voids. In Section~\ref{ss.simulation} we will show some 
simulation  examples which illustrate this statement.

\subsection{$dcx$ Ordering and Network Capacity}
\label{ss.NetCapacity}
We focus now on capacity quantifiers in the sense of communication theory.
Inequality~(\ref{e.shot-noise-LT}) with $s>0$ and 
$\Phi_1=\Phi_I$ and Poisson $\Phi_2=\Phi_\mu$
 will be used in the proof of Proposition~\ref{thm:sinr_poisson_perc} 
(see Appendix) to show that the {\em lower  level-sets
   $\{x:I_{\Phi_I}(x)\le M\}$ of the interference
field  generated by $\Phi_I$ percolate through Peierls argument}  
(cf. \cite[Proposition  14.1.4]{FnT1}) for sufficiently large $M$. 
Similarly~(\ref{e.shot-noise-LT})
with $s<0$ can be used to prove  that the {\em upper  level-sets of
  the interference 
field $I_{\Phi_I}(x)$ generated by $\Phi_I$ percolate through Peierls
argument for sufficiently large level values $M$}.

Having observed this {\em double} impact of sub-Poisson assumption 
on $\Phi_I$ it is not evident whether the threshold value $\gamma_c$
of the interference reduction factor is larger  in sub-Poisson
network than in the corresponding Poisson one.
Note that $\gamma_c$ can be related to the information-theoretic
capacity (throughput) that can be sustained on the links of the SINR graph.
In what follows we will try to explain  how sub-Poisson assumption
 can impact  some quantifiers of the network capacity.

\subsubsection{Ordering of Independent Interference Field}
  
It is quite natural to consider a network capacity characteristic
$C=f(I(x_1),\ldots, I(x_n))$ that depends on some interference
field $I$ through some $ddcx$ function $f$. We will
give a few simple examples of such characteristics in what follows.
Then, a {\em larger in $ddcx$ order interference field (more variable!)
$I(x)$ leads to larger average capacity $\E[C]$}. In particular
in the case of shot-noise interference field $I_{\Phi_I}(x)$,
by~(\ref{e.shot-noise-dcx}) we conclude that {\em larger in  $dcx$
  order pp $\Phi_I$ (clusters more!) leads to larger average
  capacity $\E[C]$}.

Let us illustrate this somewhat surprising observation by two simple
examples.
\begin{ex}[Shannon Capacity] 
\label{ex.shannon}
Let $C(I)=\log(1+F_0/(N+I))$ for some constants $F_0,N>0$ and random $I$.
Clearly this is a decreasing convex function of $I$ and
larger in convex order $I$ gives larger mean capacity $\E[C(I)]$.
The smallest value of $\E[C(I)]$ given $\E[I]$ is attained for
constant $I\equiv\E[I]$.   
\end{ex}

\begin{ex}[Outage capacity of a channel with fading]
Assume a random channel fading $F$ with convex tail
distribution function $G(t)=\Pro\{\,F> t\,\}$.
Consider $C(I)=\Pro\{\,F/(N+I)>T\,\}$ for some constant $N,T>0$.
Assuming independence of $F$ and $I$
we have $\E[C]=\E[G(T(N+I))]$, which is by our assumption expectation
of a  convex
function of $I$ and the same conclusion can be made as in
Example~\ref{ex.shannon}. The above general form of 
the expression for the outage capacity can be found in many more
detailed models based on pp, in particular in the Bipolar model of
spatial Aloha in~\cite[Chapter~16]{FnT2}.
See also~\cite{snorder} for a multidimensional version of this observation.
Similar conclusion can be made for the ergodic Shannon capacity 
$\E[\log(1+F/(N+I))]$  of 
this channel. 
\end{ex}

\subsubsection{Ordering of the Back-Bone}
Let us take one step further and consider the interference created by
the original pattern of network points (the back-bone process
$\Phi_B$) rather than an external (independent) interference field. 
In this case the interference $I_{\Phi_B}(x)$ at receiver $x$ usually has to be
considered under the so-called Palm distribution 
of $\Phi_B$ {\em given} the location of the emitter (in $\Phi_B$) of
this receiver $x$ (again see e.g. the Bipolar model of
spatial Aloha in~\cite[Chapter~16]{FnT2} for a detailed example).
The problem is that $dcx$ ordering of $\Phi_B$ implies only $idcx$
ordering of the respective Palm versions of $\Phi_B$;
cf.~\cite{snorder}. The fundamental
reason for the required ``extra'' {\em increasing} property  of the
comparable functions is that a smaller in $dcx$ pp $\Phi_B$ 
(having less clustering or even some point ``repulsion'') will  
have under Palm probability 
potentially {\em fewer points in vicinity of the conditioned point}.
This potentially decreases interference created locally near this
point, thus potentially increases our capacity characteristic.
More formally: having $idcx$ ordered of $\Phi_B$ under Palm
probability $\Pro^0$
and the capacity expressed as a $ddcx$ function of $I_{\Phi_B}$ we
cannot conclude any inequality for $\E^0[C(I_{\Phi_B})]$.

The situation is naturally inverted when we are dealing with pp which
are $dcx$ larger than Poisson pp.
A detailed analysis in~\cite{Ganti08}
of the outage capacity in the Bipolar model
generated by some  Poisson-Poisson cluster pp known as Neyman-Scott pp
(which is $dcx$ larger than Poisson pp)
confirms the above observations. Namely, for smaller transmission
distance the negative impact of clustering (locally more interferers) 
decreases the outage capacity, while for larger transmission distance 
the positive impact of interferers being more clustered increases this
capacity.

\subsubsection{Multi-hop Capacity Models} 
Percolation models have been also shown useful  to  study the 
transport (multi-hop) 
capacity of ad-hoc networks. For example, in~\cite{Franc_etal07},
by using a specific multi-hop transmission strategy that involves
percolation theory models, it was shown that a Poisson SINR network 
can achieve  capacity rate of the order of $1/\sqrt{n}$ bits
per unit of time and per node, thus closing the gap with respect to
the rate $1/\sqrt{n\log n}$ shown achievable in Poisson
networks in~\cite{Gupta00}. An interesting and open question (particularly in
view of what was shown above) is whether these capacity results can be
extended to sub-Poisson networks.

\section{Examples of Sub-Poisson Point Processes}
\label{s.Examples}

From~\cite[Section~5.2, 5.3]{snorder}, we have a rich source of examples of
Cox (doubly stochastic Poisson) pp's
comparable in $dcx$ and $ddcx$ order. 
In particular, we know that the so called L\'evy-based Cox pp 
(with Poisson-Poisson cluster pp as a special case) 
is $dcx$ larger than the Poisson pp of the same mean intensity.
Thus, they can be called as {\em $dcx$-super-Poisson}. 
However, note that any Cox pp, whose (random) realizations
of the intensity are almost surely bounded by some constant, 
can by coupled with (constructed as a subset of) 
a Poisson pp with intensity equal to
this constant. Consequently such Cox processes are trivially 
$idcx$-sub-Poisson.

\subsection{Sub-Poisson Point processes}
\label{ss.examples}
In the remaining part of this section we concentrate on 
the construction of examples of  $dcx$-sub-Poisson pp's.
\begin{ex}[Perturbed lattice]
Consider some lattice (e.g. the planar hexagonal one, 
i.e., the  usual  ``honeycomb'' model
often considered in cellular network context).
Let us ``perturb'' this ``ideal'' pattern of points as follows.
For each point of this lattice, say $z_i$, let us 
generate independently, from some given distribution, a random number,
say $N_i$, of nodes.  Moreover, instead of putting these nodes at
$X_i$, let us translate each of these nodes  
independently from $z_i$, by  vectors,  having 
some given spatial distribution (say for simplicity, of bounded support).
The resulting network, called {\em perturbed lattice},
can be seen as {\em replicating and
dispersing points from the original lattice}.
Interesting observations are as follows.
\begin{itemize}
\item
If the number of replicas $N_i$ are Poisson random variables, then the
perturbed lattice is a Poisson p.p. 
\item If moreover  the node displacement is uniform 
in the Voronoi cell of the original lattice, then the resulting
perturbed  lattice is  homogeneous Poisson
process $\Phi_\lambda$ with $\lambda$ related to $\E[N]$ and the
original lattice  density.
\item Now, if $N_i$ are convexly
smaller than some Poisson variable~\footnote{i.e., expectations of all convex
functions of $N=N_i$ are smaller than the respective expectations for
Poisson variable of the same mean; a
 special case is when $N_i=const$ is constant} 
then  the perturbed lattice is
$dcx$-sub-Poisson. This is so because the perturbed lattice pp is an
independent sum of countably many pp formed by the perturbations of
every vertex of the lattice. This means that $\Phi = \bigcup_{z_i} \Phi_{z_i}$ 
where $\Phi_{z_i}$ are independent pp. From \cite[Proposition. 3.2]{snorder}, 
we have that if $\Phi_{z_i}$ are $dcx$ ordered for every $z_i$, 
then so is $\Phi$. 
Since each of these pp is formed by $N_i$ i.i.d. perturbations, it can be shown that these pp's are $dcx$ ordered
if their respective $N_i$'s are convexly ordered. 
This follows from proving that the following function $g$ is 
convex in its argument $n$: $g(n) = \sE(f(\Phi_{z_i}(A_1),\ldots,\Phi_{z_i}(A_n))| N = n)$ for any $dcx$ function $f$, disjoint bounded Borel subsets $A_1\ldots,A_n$ and
any fixed $z_i$ of the original lattice. 
\item  Consequently, if $N_i$'s are convexly smaller than some 
Poisson random variable
and  moreover the node displacement is uniform 
in the Voronoi cell of the original lattice, then the resulting
perturbed  lattice is $dcx$-smaller than the respective homogeneous Poisson
process $\Phi_\lambda$.
\end{itemize}
\end{ex}
The interest in the above perturbed lattice models
in the networking context 
stems from the fact that they provide some spectrum of models 
that ranges from periodic grids, usually considered in cellular network
context, to ad-hoc networks almost exclusively considered under Poisson
assumptions. In Section~\ref{ss.simulation} we will show some samples
of perturbed lattices.

On the theoretical side, the interest in  perturbed lattices 
stems from their relations to
zeros of Gaussian analytic functions (GAFs)
(see~\cite{Peres05,Sodin04}).
 More precisely~\cite{Sodin06} shows that 
zeros of some  GAFs have the same distribution as points of some
 ``non-independently'' perturbed lattice.

Another class of pp's, which can be shown as {\em weakly sub-Poisson}, 
are stationary determinantal pp's. For a quick introduction
refer~\cite{Ben06}. 
\begin{ex}[Determinantal pp's] 
These pp's are 
defined by their joint intensities satisfying the following relation for all $k \geq 1$ :
$$ \rho^{(k)}(x_1,\ldots,x_k) = \det \Big (K(x_i-x_j) \Big )_{i,j}, $$
for some Hermitian, non-negative definite, locally square integrable
kernel $K : \bbR^d \to \mathbb{C}$. Then by Hadamard's inequality%
~\footnote{
$\det( a_{ij})_{ij} \leq \prod_i a_{ii}$ for Hermitian,
non-negative definite matrices $(a_{ij})$}, 
we have that {\em stationary determinantal
pp are weakly sub-Poisson}. These pp are considered as examples of pp
whose points ``exhibit repulsion''. 
Though zeros of GAFs are related to determinantal pp's, curiously enough only
zeros of i.i.d. Gaussian power series (i.e, $f(z) =
\sum_{n=0}^{\infty}a_n z^n$ , with $a_n$ i.i.d. standard complex
normal) are proved to be  determinantal pp's (see~\cite{Peres05}).    
\end{ex}

\subsection{Simulations}
\label{ss.simulation}
\begin{figure*}[!t]
\begin{center}
\begin{minipage}[b]{0.16\linewidth}
\includegraphics[width=1.08\linewidth]{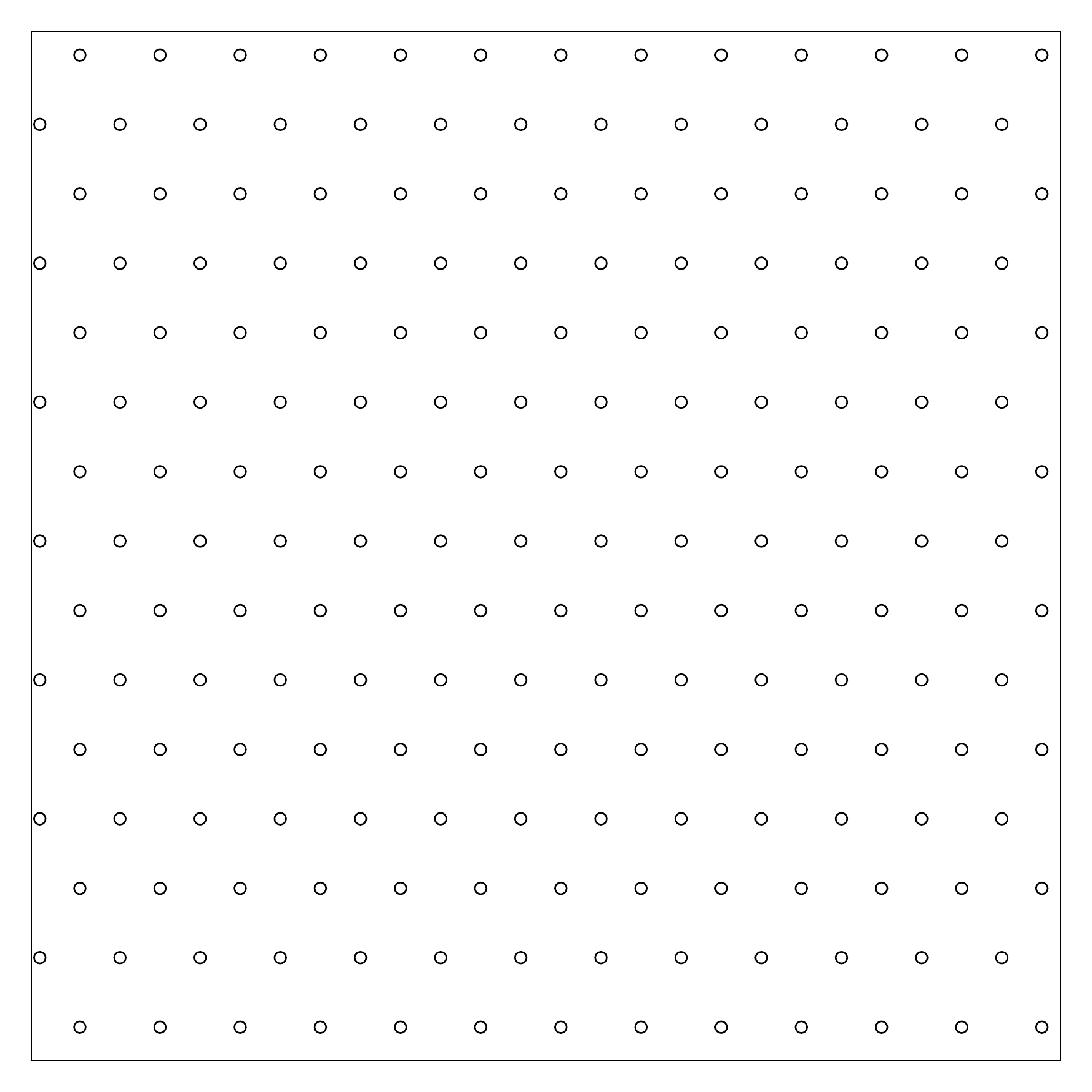}\\
\centerline{hexagonal lattice}
\end{minipage}
\begin{minipage}[b]{0.16\linewidth}
\includegraphics[width=1.08\linewidth]{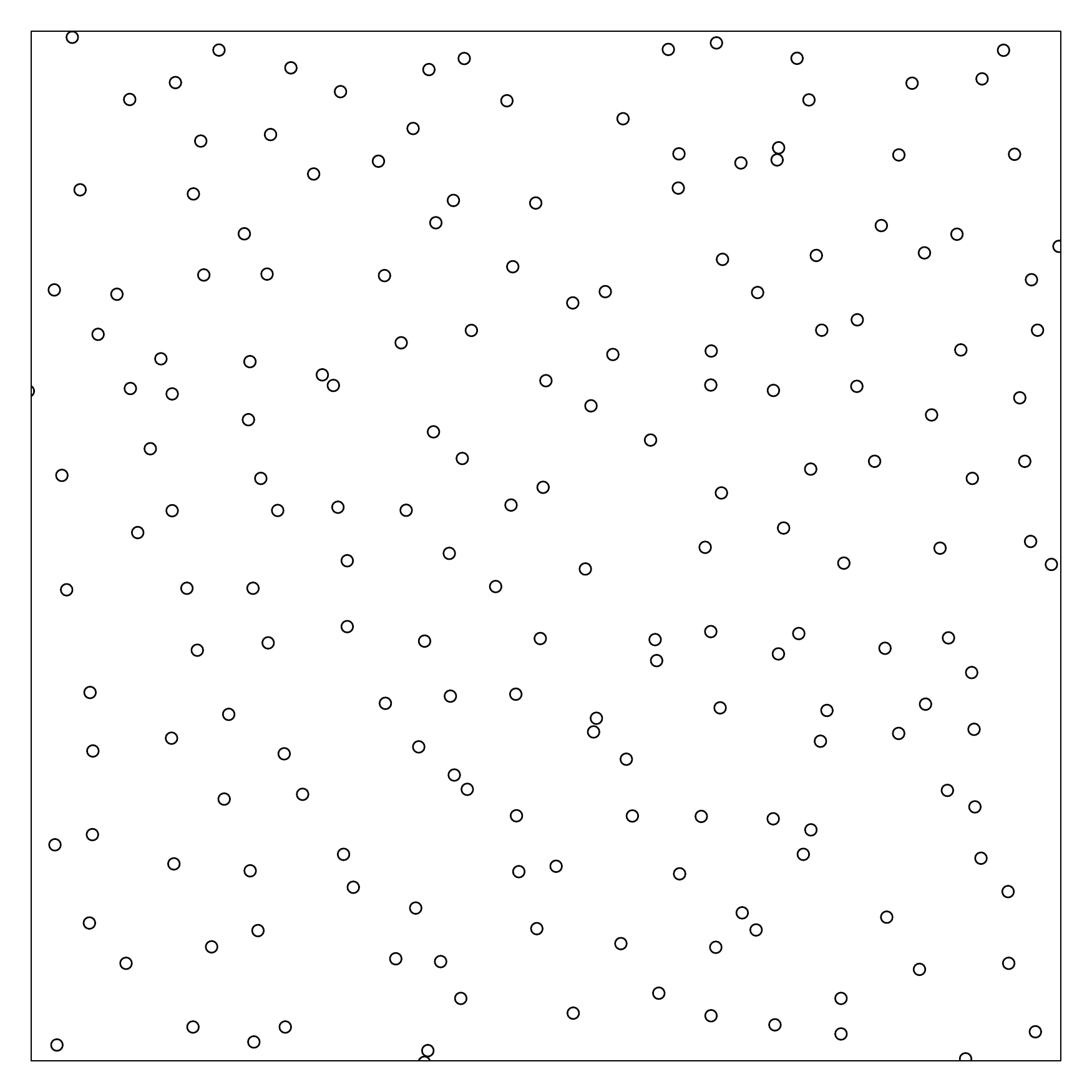}\\
\centerline{$Bin(1/1)$}
\end{minipage}
\begin{minipage}[b]{0.16\linewidth}
\includegraphics[width=1.08\linewidth]{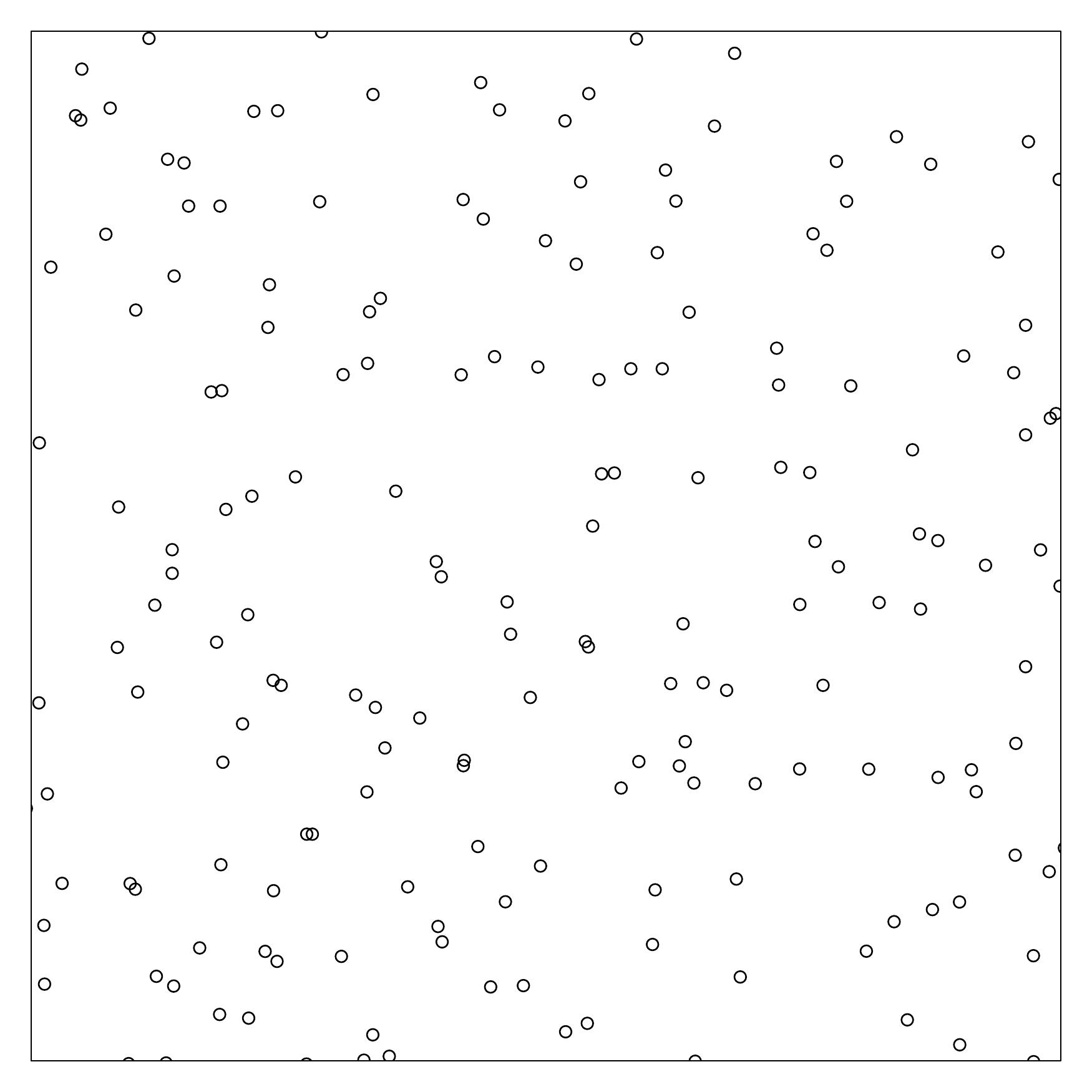}\\
\centerline{$Bin(2,1/2)$}
\end{minipage}
\begin{minipage}[b]{0.16\linewidth}
\includegraphics[width=1.08\linewidth]{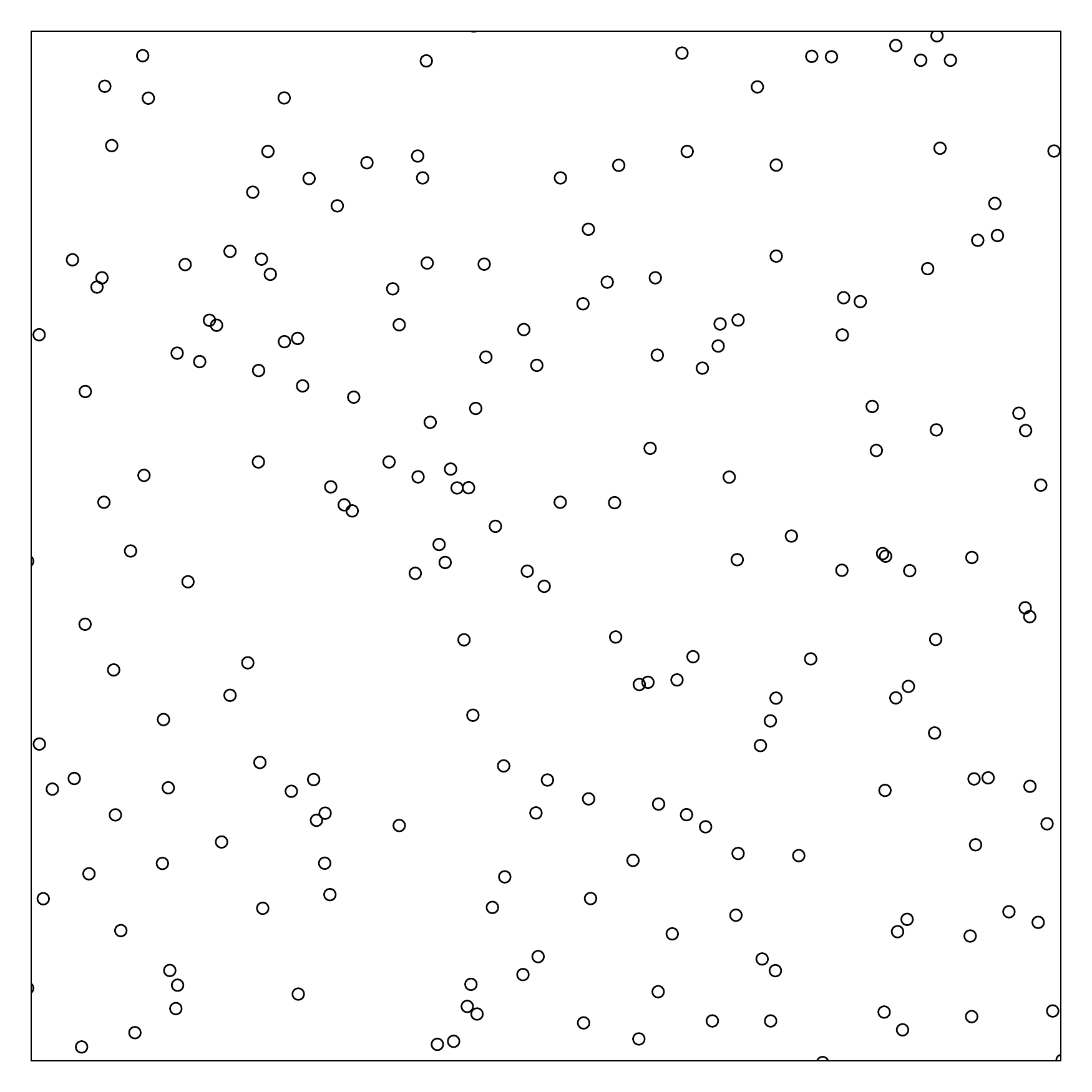}\\
\centerline{$Bin(3,1/3)$}
\end{minipage}
\begin{minipage}[b]{0.16\linewidth}
\includegraphics[width=1.08\linewidth]{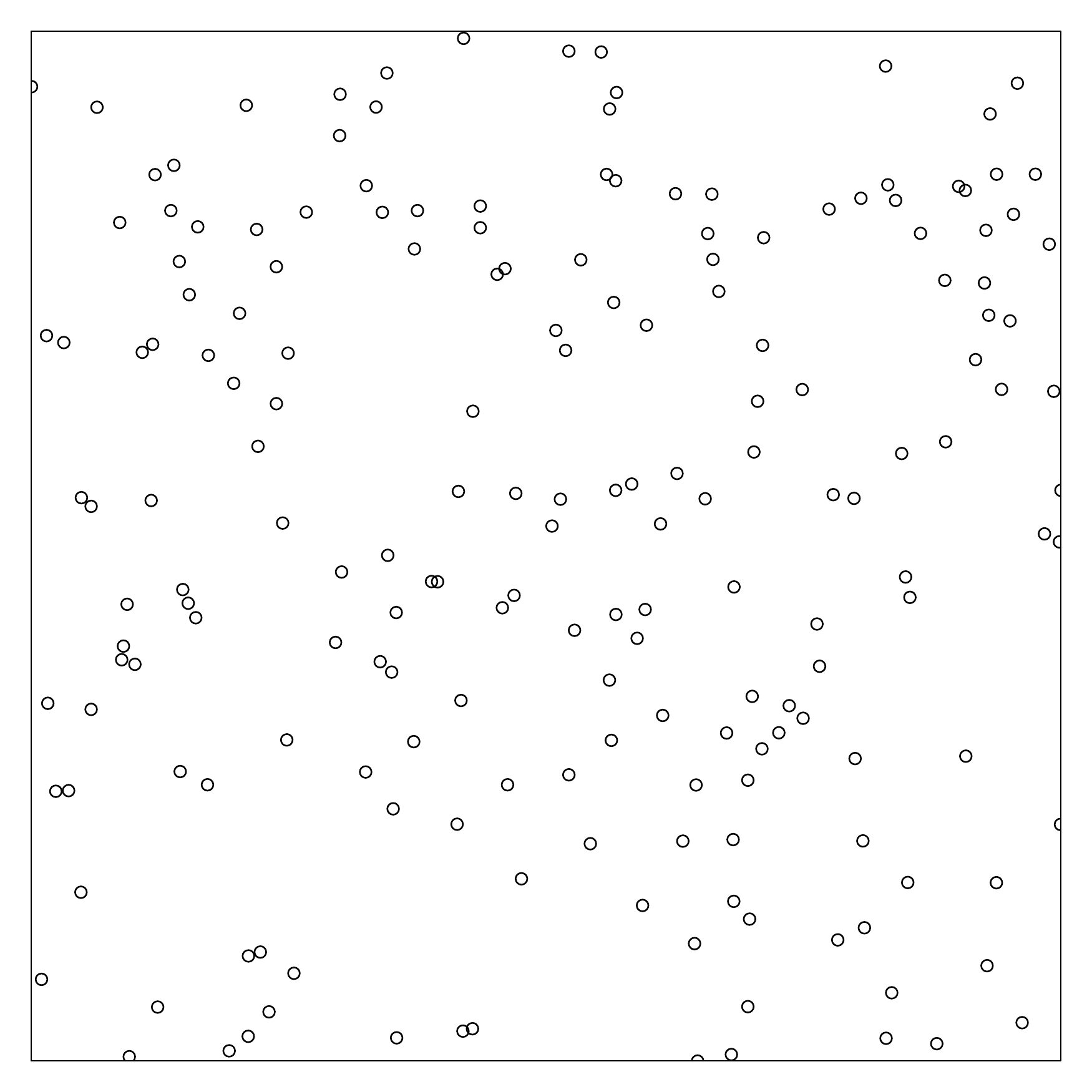}\\
\centerline{$Poi(1)$}
\end{minipage}
\begin{minipage}[b]{0.16\linewidth}
\includegraphics[width=1.08\linewidth]{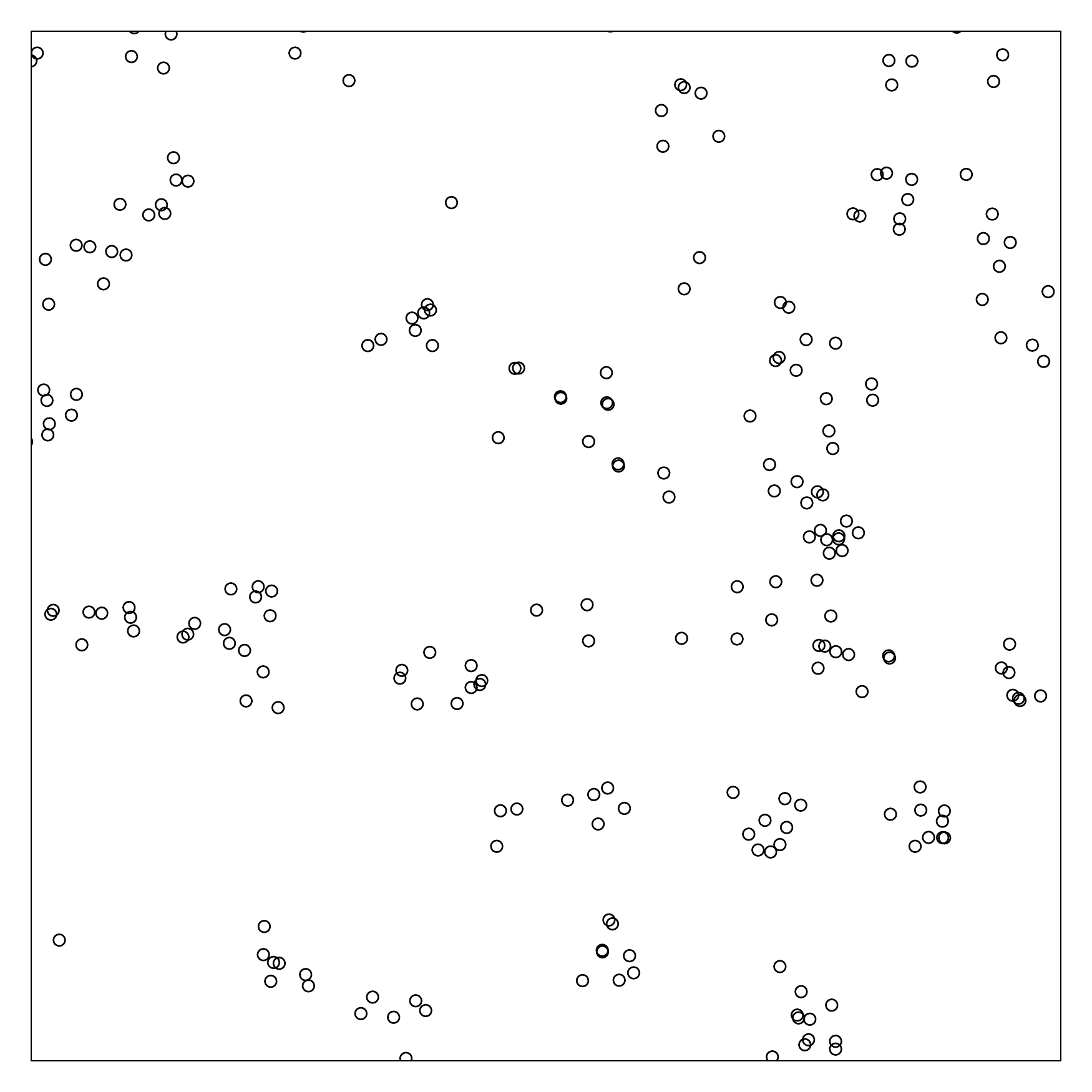}\\
\centerline{$Cox(5\times Bin(1,1/5))$}
\end{minipage}
\end{center}
\vspace{-3ex}
\caption{\label{f.Lattice1}''Unperturbed'' hexagonal lattice and 
sub-Poisson  perturbed lattices with the number of replicas $N$ having binomial distribution $B(n,1/n)$. The last figure is a super-Poisson perturbed lattice with $N$ having double stochastic Poisson (Cox) distribution of random mean $L$ having Bernoulli distribution $5 \times Bin(1,1/5)$ (i.e., $\Pr\{\,L=5\,\}=1-\Pr\{\,L=0\,\}=1/5$).}  
%
%
\begin{center}
\hbox{
\hbox{}\hspace{0.01\linewidth}
\begin{minipage}[b]{0.22\linewidth}
\includegraphics[width=1\linewidth]{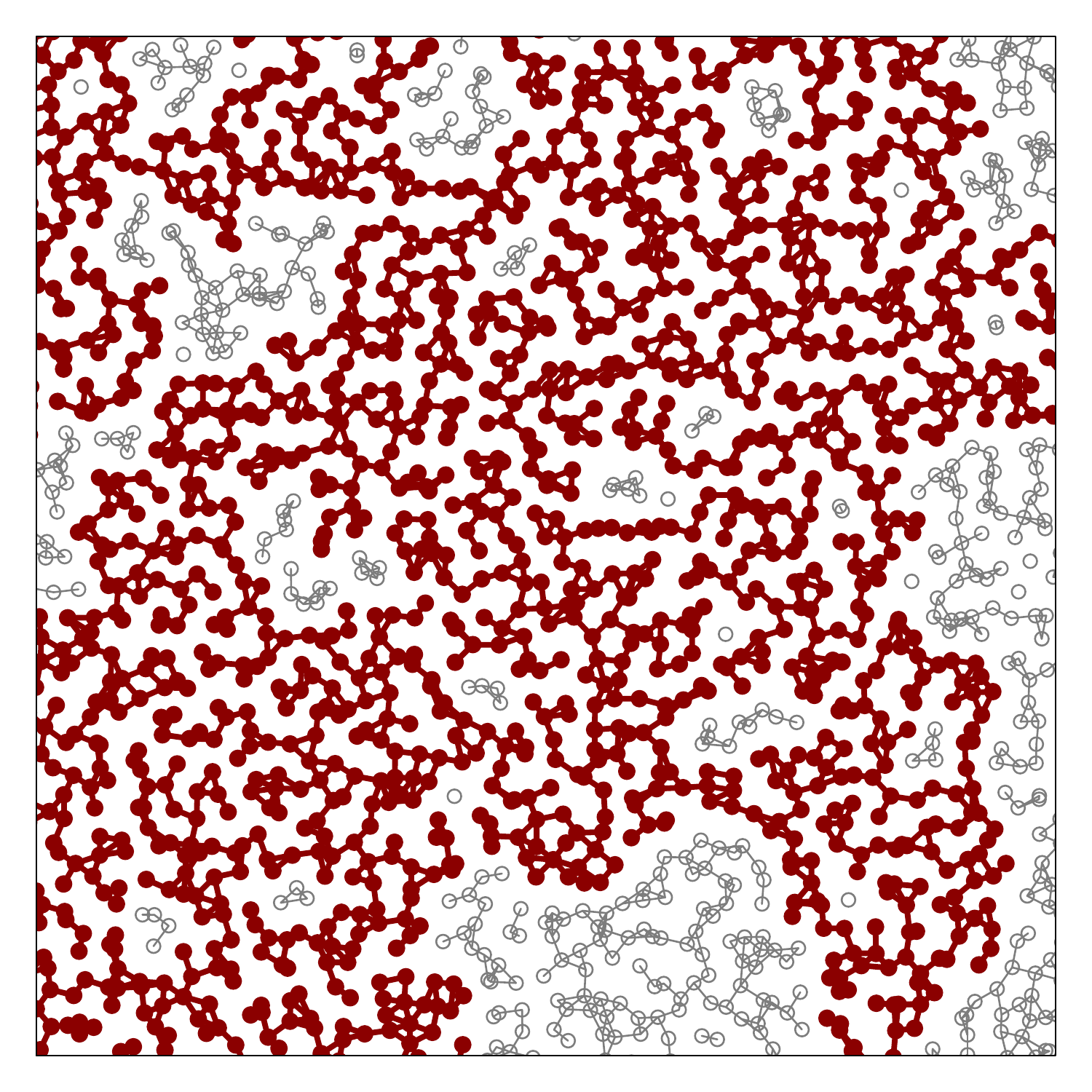}\\[-0.6\linewidth]
\hbox{}\hspace{-0.16\linewidth}
\includegraphics[width=0.70\linewidth]{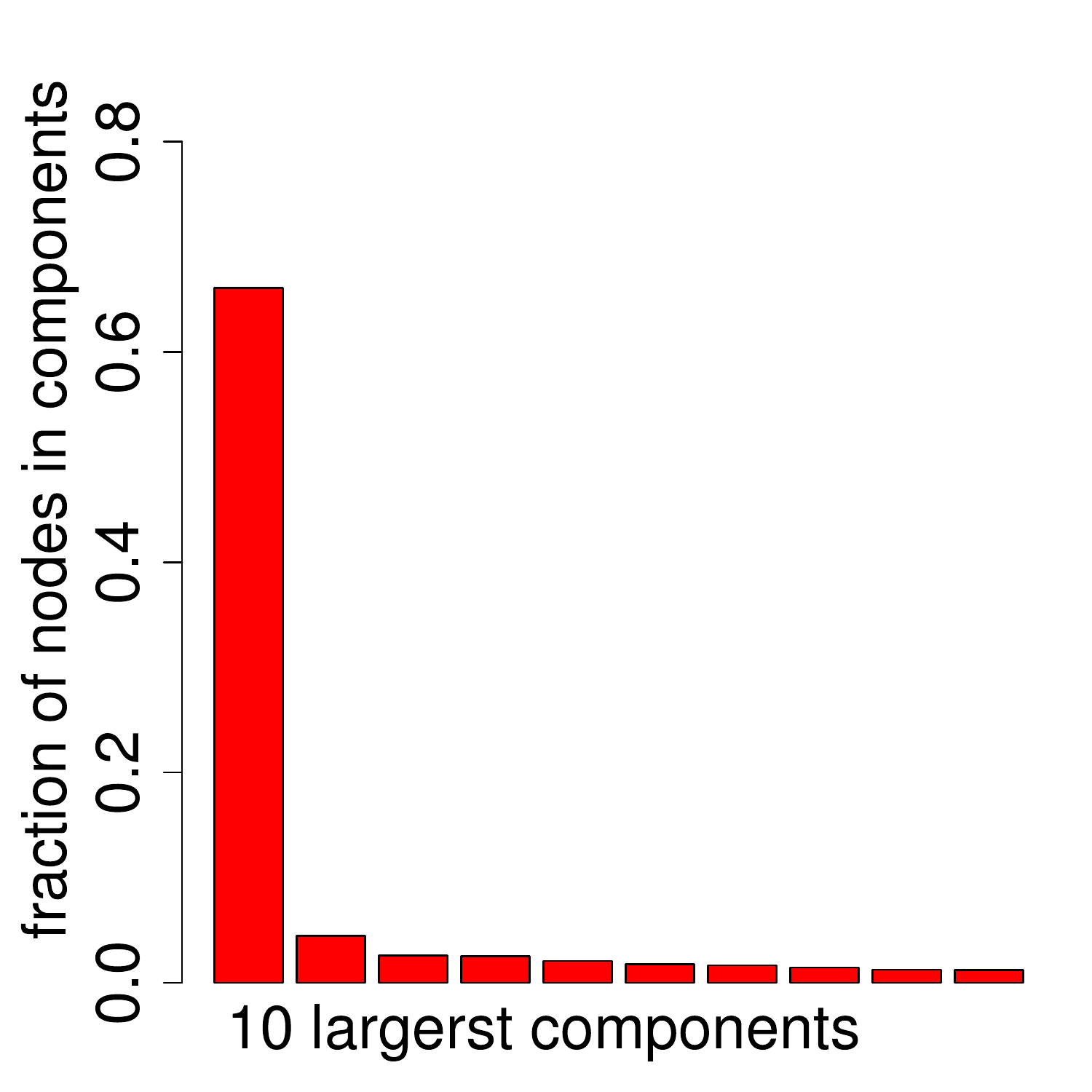}\\
\centerline{$Bin(1,1)$, $\rho=1.04$}
\end{minipage}
\hspace{0.02\linewidth}
\begin{minipage}[b]{0.22\linewidth}
\includegraphics[width=1\linewidth]{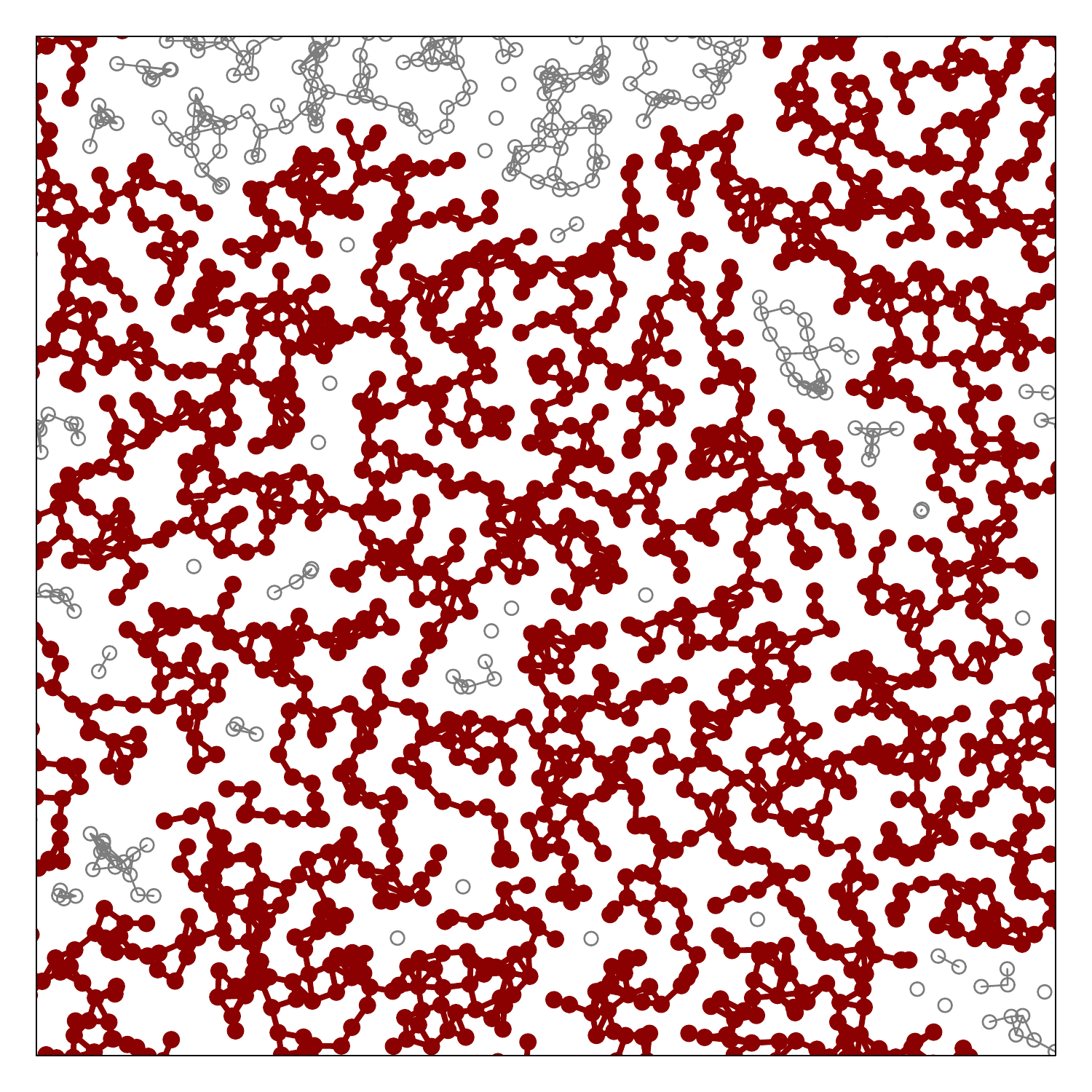}\\[-0.6\linewidth]
\hbox{}\hspace{-0.16\linewidth}
\includegraphics[width=0.70\linewidth]{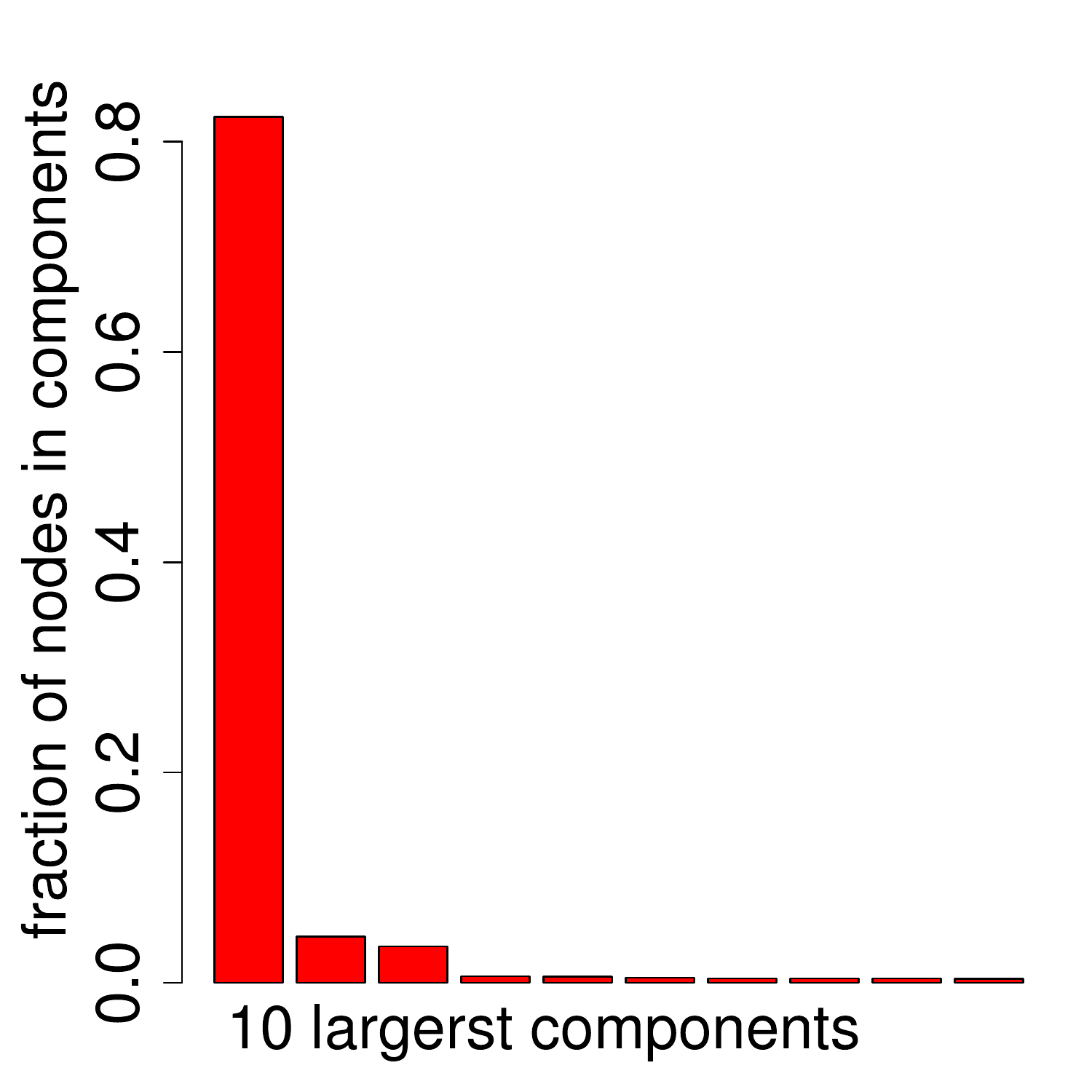}\\
\centerline{$Bin(2,1/2)$, $\rho=1.07$}
\end{minipage}
\hspace{0.02\linewidth}
\begin{minipage}[b]{0.22\linewidth}
\includegraphics[width=1\linewidth]{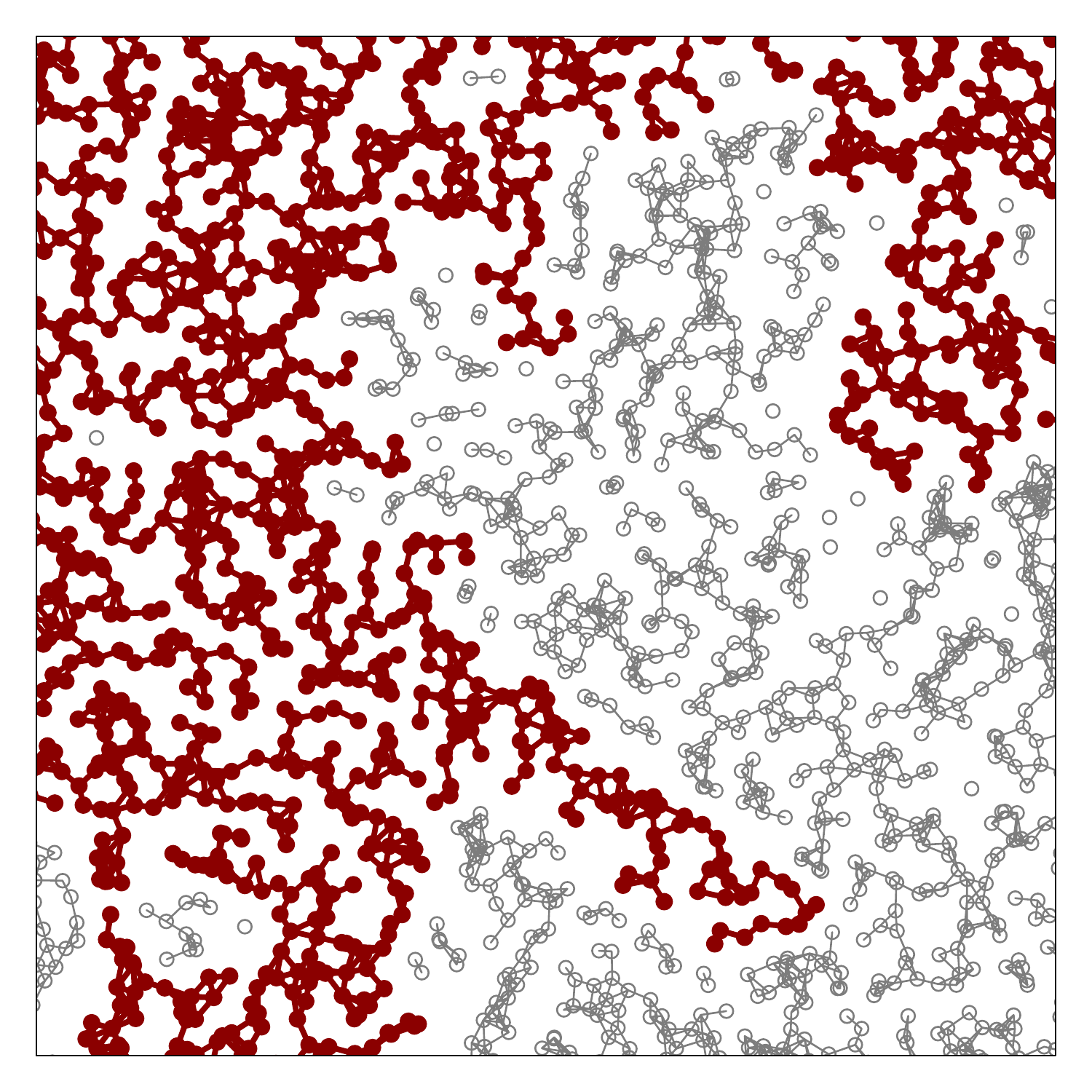}\\[-0.6\linewidth]
\hbox{}\hspace{-0.16\linewidth}
\includegraphics[width=0.70\linewidth]{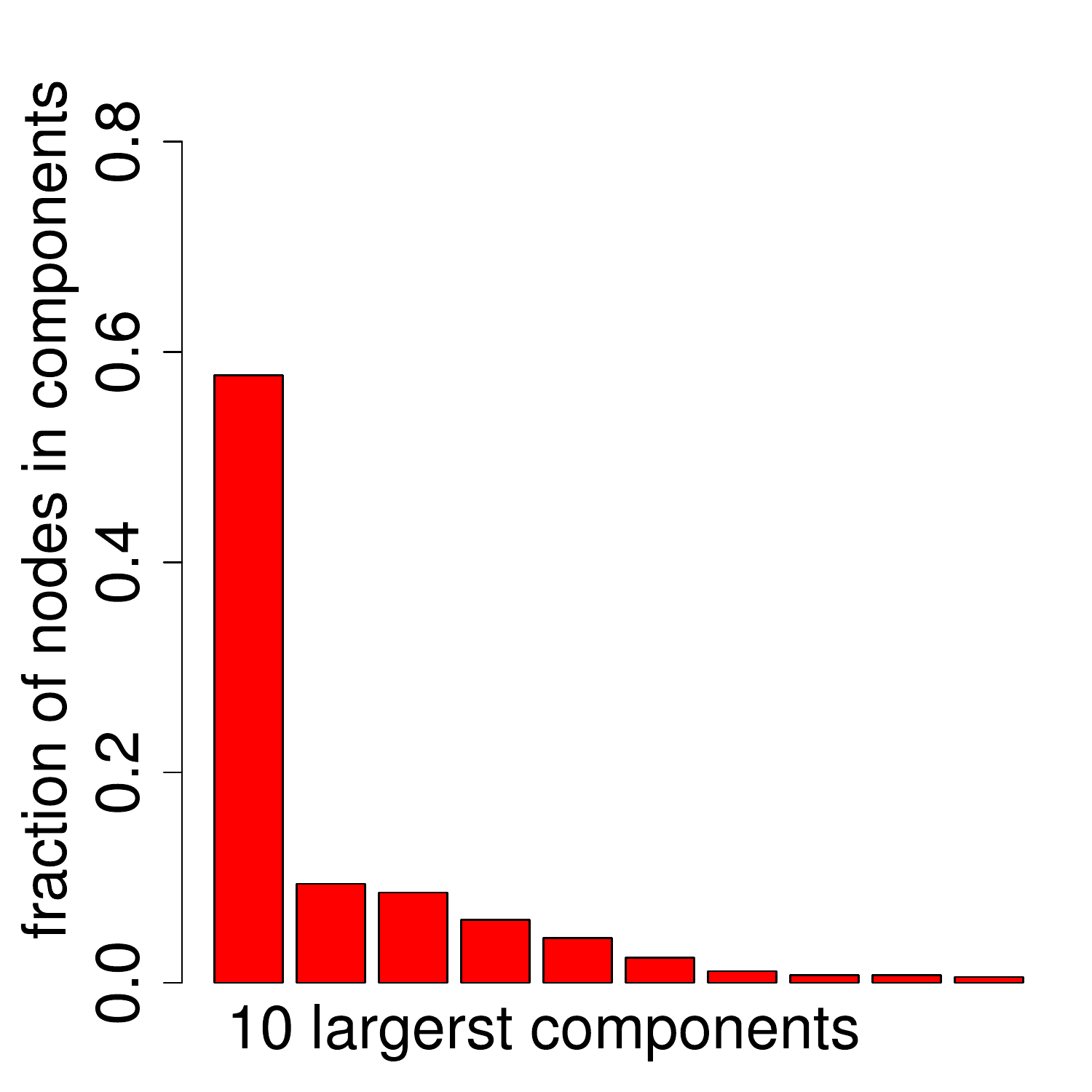}\\
\centerline{$Bin(3,3)$, $\rho=1.09$}
\end{minipage}
\hspace{0.02\linewidth}
\begin{minipage}[b]{0.22\linewidth}
\includegraphics[width=1\linewidth]{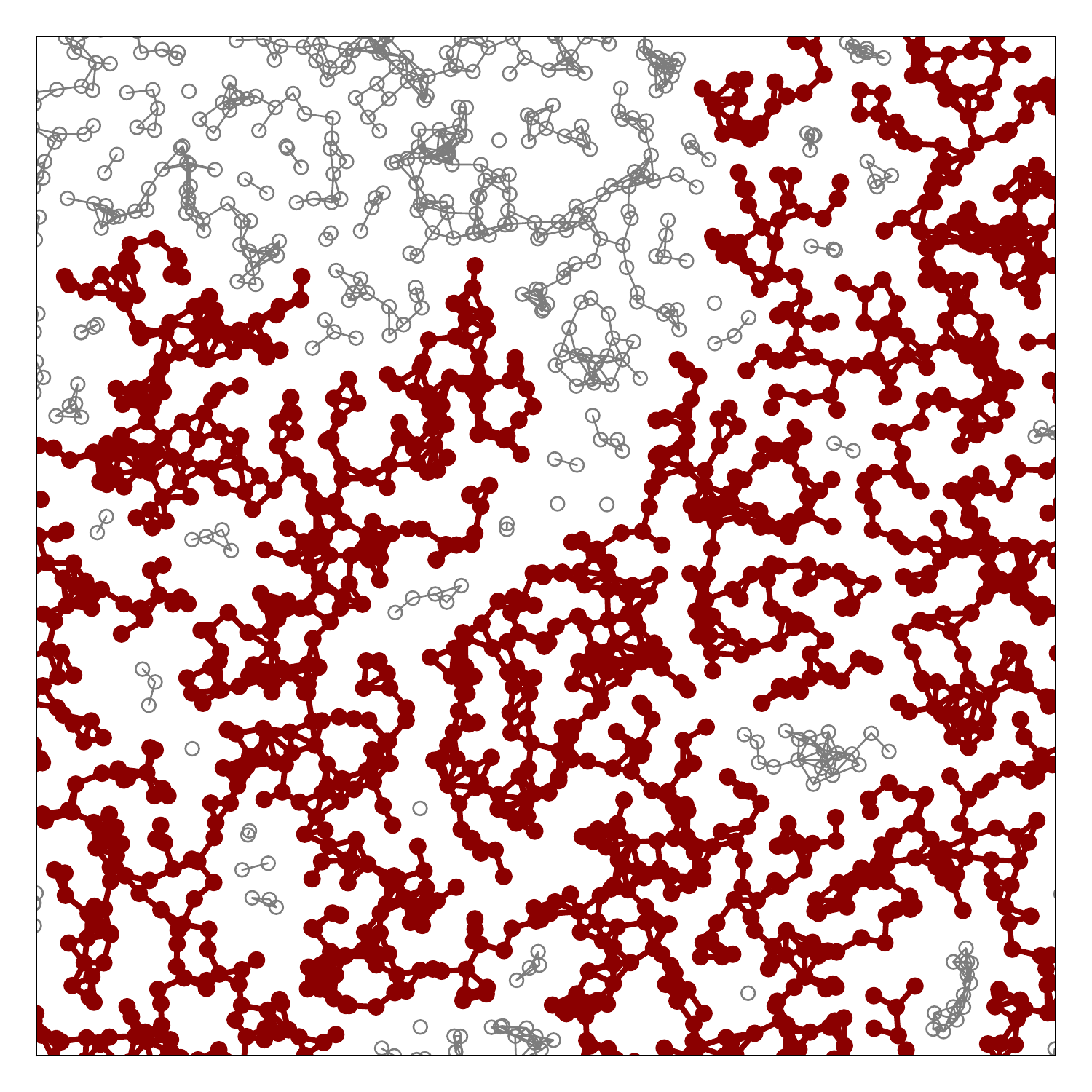}\\[-0.6\linewidth]
\hbox{}\hspace{-0.16\linewidth}
\includegraphics[width=0.70\linewidth]{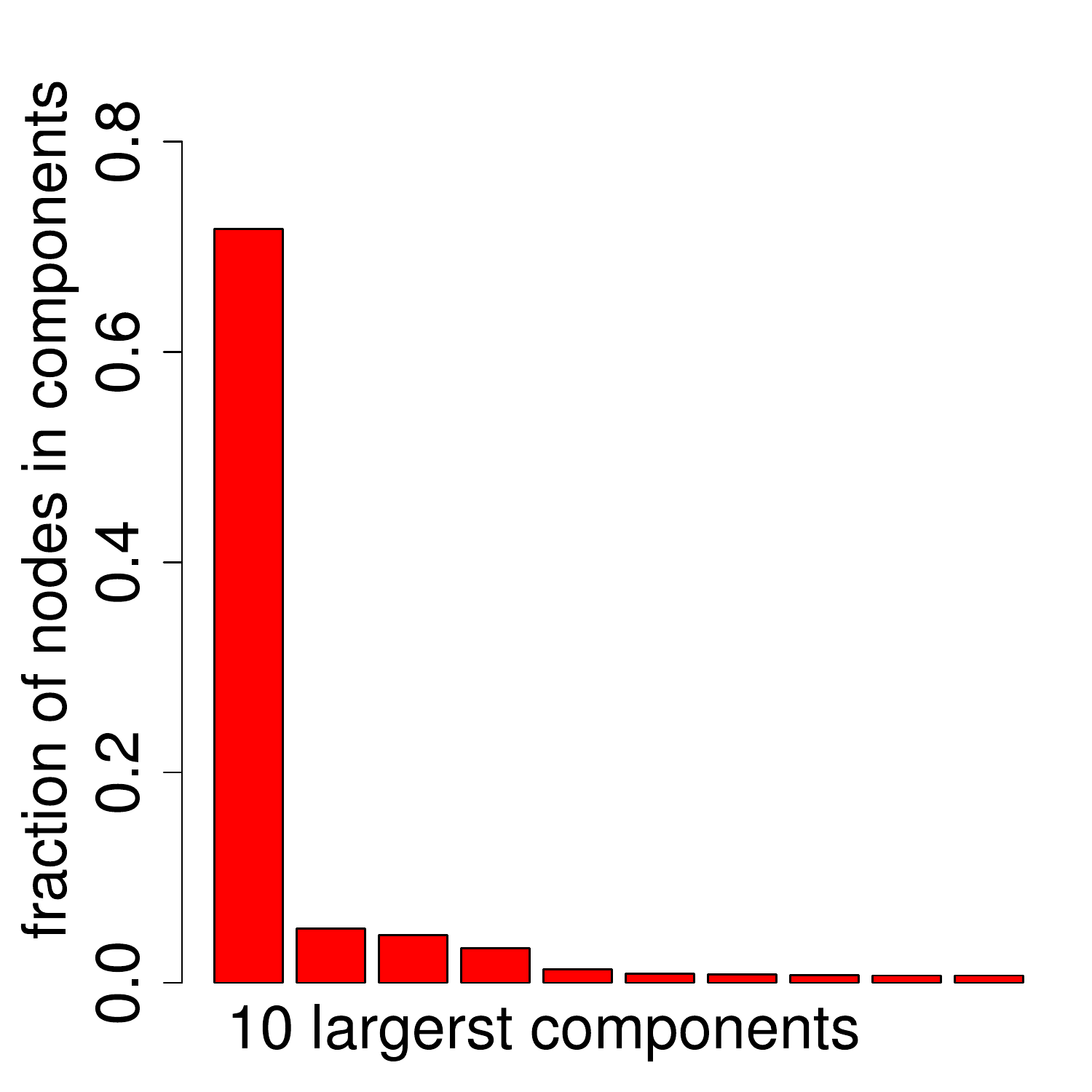}\\
\centerline{$Poi(1)$, $\rho=1.12$}
\end{minipage}}
\end{center}
\vspace{-5ex}
\caption{\label{f.Gilbert1}
Gilbert graph with communication range $\rho$
and nodes form a perturbed lattice pp with Binomial $Bin(n,1/n)$
number of replicas uniformly distributed in hexagonal cells. The
largest component in the simulation window is highlighted. Bar-plots
show the fraction of nodes in ten largest components.}
\vspace{-3ex}
\end{figure*}

Now, we will show some examples of simulated 
patterns of perturbed lattices. We consider hexagonal pattern of 
original (unperturbed) points (cf. Figure~\ref{f.Lattice1}, upper-left plot).
The replicas are always displaced uniformly in the Voronoi cell 
of the given point of the original point. Consequently, when the
numbers of replicas $N$ have Poisson distribution then the corresponding
perturbed lattice is Poisson pp. (cf. Figure~\ref{f.Lattice1}, lower-middle
plot). We always take $\E[N]=1$.
In order to generate sub-Poisson pp's we take $N$ having Binomial
distribution $Bin(n,1/n)$. This family of distributions can be shown to be {\em convexly increasing in
$n$} and {\em convex upper-bounded by Poisson} $Poi(1)$ distribution.
Note on Figure~\ref{f.Lattice1}
that in the case $Bin(1,1)$, where we have only
point displacement (no replications), the resulting perturbed lattice
looks much more ``regular'' (less visible voids, less clusters)
than Poisson pp. This  regularity
becomes less evident when $n$ increases and already for $n\ge3$ 
it is difficult to distinguish the perturbed lattice patterns from
Poisson pp. However, a more detailed study  of 
Gilbert graphs on them (not presented
here) up to approximately $n=5$ show that they
still have significantly different clustering structure than that of
Poisson pp.

The lower-right plot on Figure~\ref{f.Lattice1} shows an example of
perturbed lattice $dcx$ {\em larger} than Poisson pp. 
It is a  doubly stochastic Poisson pp (Cox pp), where the mean value $L$ 
of the number of replicas $N$ is first sampled independently for each 
original lattice point, with 
$\Pr\{\,L=a\,\}=1-\Pr\{\,L=0\,\}=1/a$, for some $a\ge1$  
and then, given $L$, $N$ has
Poisson distribution $Poi(L)$. This construction is an example of
the Ising-Poisson cluster pp considered in~\cite[Section~5.1]{snorder}.
Another possibility to generate the  number of replicas $N$ convexly
larger than Poisson variable is to ``scale-up''  the Poisson variable 
of a smaller intensity, i.e., to take 
$N=nN'$ where $N'\sim Poi(1/n)$ for arbitrary $n\ge1$.

In Figure~\ref{f.Gilbert1} we consider a larger
simulation window comprising  about $40^2$ points of the original
hexagonal lattice. We generate  points of the perturbed lattices
with $N\sim Bin(n,1/n)$ for $n=1,\ldots,3$ and $n=20$ to
approximate Poisson pp. For each $n$ we show
the  Gilbert graph with the communication radius 
 $\rho$ for which the largest component  (the highlighted one)
starts to significantly out-number all other components.
More precisely for each given simulated pattern of points we find 
the smallest $\rho=\rho(n)$ (up to the second decimal place) 
such that the largest component in the simulation window 
contains about 60\% of the simulated points 
(cf. the bar-plots showing the empirical
fraction of the number of points in 10 largest components). 
Observing values of $\rho(n)$ we conjecture that  {\em the critical
  radius $\rho_c(n)$ for the 
  percolation of the Gilbert graph on the considered perturbed
  lattices  is increasing in $n$}. 
Note that for the unperturbed lattice $\rho_c=1$ (the
distance between adjacent nodes in our hexagonal lattice) and for
Poisson pp of the same intensity  $\rho_c$ is known to be close
to~$1.112$.  A more exhaustive numerical study can be found in the
extended version of the paper~\cite{subpoisson-ext}.

\section{Concluding Remarks}
\label{s.Conclusions}
We have extended two results on nontrivial phase transition 
in  percolation of Gilbert's and SINR
graphs from Poisson to  $dcx$-sub-Poisson pp's.
This means that, regarding existence of this phase transition,
Poisson pp represents a worst-case scenario
within this class of pp's.

A natural question in this context is as follows. Consider
$\Phi_1\leq_{dcx}\Phi_2$. Does this ordering 
imply that the corresponding critical  communication ranges 
in the Gilbert's model are ordered as well $\rho_c(\Phi_1)\le\rho_c(\Phi_2)$?
The answer in the full generality might be negative, as we expect some
counterexamples of $dcx$-super-Poisson pp's with the critical
communication ranges degenerating to $0$ or $\infty$.
However, we conjecture that the critical communication 
range of a stationary $dcx$-sub-Poisson pp's of intensity $\lambda$
is smaller than that of a Poisson pp of intensity $\lambda$.

In this paper we have also discussed the impact of $dcx$ ordering
on some local quantifiers of the capacity in networks. 
An interesting open question in this context 
is the relation between $dcx$ ordering
and the transport capacity in these networks. 

Finally, we brought attention to the so called perturbed lattices,
which can provide a spectrum of pp's  monotone in $dcx$ order
that ranges from periodic grids, usually considered in cellular network
context, to Poisson pp's and to various more clustered pp's including
doubly stochastic Poisson pp's.

\appendix
\label{s.appendix} 
\renewcommand{\theequation}{A.\arabic{equation}}
\section{Appendix}
In this section we prove   our main results
of this paper.

\begin{proof}({\em of Proposition}~\ref{thm:perc_sub-Poisson_pp})
Denote by $\bbZ^d(r)$, the discrete graph formed by the vertexes of
$r\bbZ^d=\{r z:z\in \bbZ^d\}$, where $\bbZ^d$ is the  $d$-dimensional
integer lattice,  
and edges between $z_i,z_j\in\bbZ^d(r)$ such that $\{z_i + 
[-\frac{r}{2},\frac{r}{2}]^d\} \cap \{z_j +
[-\frac{r}{2},\frac{r}{2}]^d\} \neq \emptyset.$ Now we define a site
percolation on $\bbZ^d(r)$ induced by the pp $\Phi$ : $X_r(z) =
\mathbf{1} [\Phi(Q_r(z)) \geq 1]$ where $Q_r(z) = 
(-\frac{r}{2},\frac{r}{2}]^d.$ Note that when $X_{2r}(z)$ does not 
  percolate, $C(\Phi,r)$ does not percolate and if
  $X_{\frac{r}{\sqrt{d}}}(z)$ percolates, so does $C(\Phi,r)$.

The standard technique to show non-percolation is to show that the
expected number of (self-avoiding) paths of length $n$ starting at the
origin in the random sub-graph induced
by  opened sites of the percolation model tends to zero as $n\to\infty$.
The probability of a  path
$(z_1,\ldots,z_n)$ of length $n$ in $\bbZ^d(2r)$ being supported by 
open sites
is  $\Pro\{\prod_i \Phi(Q_{2r}(z_i)) \geq 1\} \leq 
\sE(\prod_i \Phi(Q_{2r}(z_i))) \leq (\lam (2r)^d)^n$ where the first inequality is
due to Markov's inequality and the last inequality is by weak
sub-Poisson property of $\Phi$, which in particular is  
implied by  $idcx$-sub-Poisson property
(cf.~(\ref{e.weakly-sub-Poi})).
Since the number of paths starting at the origin, of length
$n$ in $\bbZ^d(2r)$, is bounded by 
$(3^d-2)^n$, we have that the expected number of
paths of length $n$ is at most  $((3^d-2)\lam (2r)^d)^n$ 
and this tends to $0$ for $r$ small enough. 
This shows that there exists $c(\lam) > 0$ (depending also on the
dimension)  such
that $r_c(\Phi) \geq c(\lam)$. 

For the upper bound, we use the Peierls argument (cf.~\cite[
  pp.~17--18]{Gr99}) on the site percolation model induced by
$X_{\frac{r}{\sqrt{d}}}(z)$ on $\bbZ^d(\frac{r}{\sqrt{d}}).$ To use
this argument, one needs to estimate the probability of the
site-percolation model not intersecting a path
$(z_1,\ldots,z_n)$. This probability can be expressed as 
$\Pro\{\Phi(\bigcup_{i=1}^n Q_{\frac{r}{\sqrt{d}}}(z_i)) = 0\}$.

and is smaller than $\exp\{- \lam n (\frac{r}{\sqrt{d}})^d\}$
by the $ddcx$-sub-Poisson property of $\Phi$ and
\cite[Proposition~4.1]{snorder}. Now by choosing $r$ large enough this
probability can 
be made as small as we wish and so the Peierls argument can be used.  
\end{proof}

\begin{proof}({\em of Proposition}~\ref{thm:sinr_poisson_perc})
We follow the proof  given in~\cite{Dousse_etal}.
Assuming $\lambda>\lambda_c(\rho_l)$, one observes first
that the graph $G(\lambda)$ also percolates
with any slightly larger constant noise
$N'=N+\delta'$, for some $\delta'>0$.
Essential for the original proof of the result 
is to show that the
level-set $\{x:I_{\Phi_I}(x)\le M\}$ of the interference
field  percolates (contains an infinite connected component) for sufficiently
large $M$. Suppose that it is true. 
Then taking $\gamma=\delta'/M$ one has 
percolation of the level-set $\{y:\gamma I_{\Phi_I}(y)\le \delta'\}$.
The main difficulty consists in showing that 
$G(\lambda)$ with noise $N'=N+\delta'$ 
percolates {\em within} an infinite connected component of 
$\{y:I_{\Phi_I}(y)\le \delta'\}$. This was done in~\cite{Dousse_etal}, 
by mapping both models $G(\lambda)$ and the level-set of the
interference to a discrete lattice
and showing that both discrete approximations not only percolate 
but actually satisfies a stronger, sufficient condition for
percolation, related to the Peierls argument~\cite[Proposition 14.1.4]{FnT1}. 
We follow exactly the same steps and the only fact that we have to
prove,  regarding  the interference, is that there exists a constant
$\epsilon<1$ such that for arbitrary $n\ge1$ and 
arbitrary choice of locations $x_1,\ldots,x_n$ one has 
$\Pro\{I_{\Phi_I}(x_i)> M, \, i=1,\ldots, n\}\le \epsilon^n$.
To this regard, as in~\cite{Dousse_etal}, 
using the Chernoff bound we dominate this probability by 
$e^{-snM}\E[\exp\{\sum_{i=1}^nsI_{\Phi_I}(x_i)\}]$ with arbitrary
$s>0$. The crucial observation for our extension of the original proof
is that $f(u_1,\ldots,u_n)=\exp[s(u_1+\ldots+u_n)]$ is an
$idcx$ function. By the assumption $\Phi_I\le_{idcx}\Phi_\mu$ 
and~\cite[Theorem~2.1]{snorder} we have 
$$\E\Bigl[\exp\Bigl\{s\sum_{i=1}^nI_{\Phi_I}(x_i)\Bigr\}\Bigr]\le
\E\Bigl[\exp\Bigl\{s\sum_{i=1}^nI_{\Phi_\mu}(x_i)\Bigr\}\Bigr]$$ 
and we can use the explicit form
of the Laplace transform of the Poisson shot-noise (the
right-hand-side in the above inequality) to prove, exactly as
in~\cite{Dousse_etal}, that for sufficiently small $s$
it is not larger than $K^n$ for some
constant $K$  which depends on $\mu$ but not on $M$.
This completes the proof.
\end{proof}

\begin{proof}({\em of Proposition}~\ref{thm:perc_sinr_sub-Poisson})
In this scenario, increased power is equivalent to increased radius in
the Gilbert's model associated with the SINR model. From this
observation, it follows that by using the discrete mapping and
arguments as in the proof of 
Proposition~\ref{thm:perc_sub-Poisson_pp}, we obtain that with
increased power the associated Gilbert's model percolates. Then, we
use the approach from the proof of
Proposition~\ref{thm:sinr_poisson_perc} to obtain 
a $\gamma > 0$ such that the SINR network percolates as well. 
\end{proof}

\singlespacing
\bibliographystyle{plainnat}
{\footnotesize 
\bibliography{subpoisson}

\begin{thebibliography}{17}
\providecommand{\natexlab}[1]{#1}
\providecommand{\url}[1]{\texttt{#1}}
\expandafter\ifx\csname urlstyle\endcsname\relax
  \providecommand{\doi}[1]{doi: #1}\else
  \providecommand{\doi}{doi: \begingroup \urlstyle{rm}\Url}\fi

\bibitem[Baccelli and B{\l}aszczyszyn(2009{\natexlab{a}})]{FnT1}
F.~Baccelli and B.~B{\l}aszczyszyn.
\newblock \emph{Stochastic Geometry and Wireless Networks, Volume I ---
  Theory}, volume 3, No 3--4 of \emph{Foundations and Trends in Networking}.
\newblock NoW Publishers, 2009{\natexlab{a}}.

\bibitem[Baccelli and B{\l}aszczyszyn(2009{\natexlab{b}})]{FnT2}
F.~Baccelli and B.~B{\l}aszczyszyn.
\newblock \emph{Stochastic Geometry and Wireless Networks, Volume II ---
  Applications}, volume 4, No 1--2 of \emph{Foundations and Trends in
  Networking}.
\newblock NoW Publishers, 2009{\natexlab{b}}.

\bibitem[Ben~Hough et~al.(2006)Ben~Hough, Krishnapur, Peres, and Virag]{Ben06}
J.~Ben~Hough, M.~Krishnapur, Y.~Peres, and B.~Virag.
\newblock Determinantal processes and independence.
\newblock \emph{Probability Surveys}, 3:\penalty0 206--229, 2006.

\bibitem[B{\l}aszczyszyn and Yogeshwaran(2009)]{snorder}
B.~B{\l}aszczyszyn and D.~Yogeshwaran.
\newblock Directionally convex ordering of random measures, shot-noise fields
  and some applications to wireless networks.
\newblock \emph{Adv. Appl. Probab.}, 41:\penalty0 623--646, 2009.

\bibitem[B{\l}aszczyszyn and Yogeshwaran(2010)]{subpoisson-ext}
B.~B{\l}aszczyszyn and D.~Yogeshwaran.
\newblock Connectivity in sub-{P}oisson networks.
\newblock http://hal.inria.fr/inria-00497707, 2010.

\bibitem[Dousse et~al.(2005)Dousse, Baccelli, and Thiran]{Dousse_etal_TON}
O~Dousse, F.~Baccelli, and P~Thiran.
\newblock Impact of interferences on connectivity in ad-hoc networks.
\newblock \emph{IEEE/ACM Trans. Networking}, 13:\penalty0 425--543, 2005.

\bibitem[Dousse et~al.(2006)Dousse, Franceschetti, Macris, Meester, and
  Thiran]{Dousse_etal}
O.~Dousse, M.~Franceschetti, N.~Macris, R.~Meester, and P.~Thiran.
\newblock Percolation in the signal to interference ratio graph.
\newblock \emph{Journal of {A}pplied {P}robability}, 43\penalty0 (2):\penalty0
  552--562, 2006.

\bibitem[Franceschetti et~al.(2007)Franceschetti, Dousse, C., and
  P.]{Franc_etal07}
M.~Franceschetti, O.~Dousse, Tse D.~N. C., and Thiran P.
\newblock Closing the gap in the capacity of wireless networks via percolation
  theory.
\newblock \emph{IEEE Trans. Inf. Theory}, 53:\penalty0 1009--1018, 2007.

\bibitem[Ganti and Haenggi(2009)]{Ganti08}
R.~Ganti and M.~Haenggi.
\newblock Interference and outage in clustered wireless ad hoc networks.
\newblock \emph{IEEE Tr. Inf. Theory}, 55:\penalty0 4067--4086, 2009.

\bibitem[Gilbert(1961)]{Gilbert61}
E.~N. Gilbert.
\newblock Random plane networks.
\newblock \emph{SIAM J.}, 9:\penalty0 533--543, 1961.

\bibitem[Grimmett(1989)]{Gr99}
G.~R. Grimmett.
\newblock \emph{Percolation}.
\newblock Springer-Verlag, New York, 1989.

\bibitem[Gupta and Kumar(2000)]{Gupta00}
P.~Gupta and P.~R. Kumar.
\newblock The capacity of wireless networks.
\newblock \emph{IEEE Transactions on Information Theory}, 46\penalty0
  (2):\penalty0 388--404, 2000.

\bibitem[Matheron(1975)]{Matheron75}
G.~Matheron.
\newblock \emph{Random Sets and Integral Geometry}.
\newblock John Wiley \& Sons, London, 1975.

\bibitem[Meester and Roy(1996)]{MR96}
R.~Meester and R.~Roy.
\newblock \emph{Continuum percolation}.
\newblock Cambridge University Press, 1996.

\bibitem[Peres and Virag(2005)]{Peres05}
Y.~Peres and B.~Virag.
\newblock Zeros of the i.i.d. {G}aussian power series: a conformally invariant
  determinantal process.
\newblock \emph{Acta Mathematica}, 194:\penalty0 1--35, 2005.

\bibitem[Sodin and Tsirelson(2004)]{Sodin04}
M.~Sodin and B.~Tsirelson.
\newblock Random complex zeroes; {I.} asymptotic normality.
\newblock \emph{Israel J. Math.}, 144:\penalty0 125--149, 2004.

\bibitem[Sodin and Tsirelson(2006)]{Sodin06}
M.~Sodin and B.~Tsirelson.
\newblock Random complex zeroes; {II.} perturbed lattice.
\newblock \emph{Israel J. Math.}, 152:\penalty0 105--124, 2006.

\end{thebibliography}
}

\end{document}